\documentclass[a4paper,12pt] {article}
\usepackage{graphics}
\usepackage{graphicx}
\usepackage{epsfig}
\usepackage{mathrsfs}
\usepackage{amssymb}
\usepackage{amsmath}
\usepackage{pstricks}
\usepackage{lmodern}
\usepackage{xcolor,rotating,epic,eepic}

\usepackage[latin1]{inputenc}
\usepackage[T1]{fontenc}

\usepackage{color}
\usepackage{pst-all}

\begin{document}

\begin{center}
\Large{\textbf{Performance, Robustness and Sensitivity Analysis of \\the Nonlinear Tuned Vibration Absorber}}\vspace{1cm}
\end{center}

\begin{center}
{T. Detroux, G. Habib, L. Masset, G. Kerschen\\\vspace{0.8cm}

\small Space Structures and Systems Laboratory\\
Department of Aerospace and Mechanical Engineering\\
University of Liège, Liège, Belgium \\
E-mail: tdetroux, giuseppe.habib, luc.masset, g.kerschen@ulg.ac.be\\\vspace{0.5cm} \vspace{1cm}

\rule{0.85\linewidth}{.3pt}
\vspace{-0.5cm}
\begin{abstract}
The nonlinear tuned vibration absorber (NLTVA) is a recently-developed nonlinear absorber which generalizes Den Hartog's equal peak method to nonlinear systems. 
If the purposeful introduction of nonlinearity can enhance system performance, it can also give rise to adverse dynamical phenomena, 
including detached resonance curves and quasiperiodic regimes of motion. Through the combination of numerical continuation of periodic solutions, 
bifurcation detection and tracking, and global analysis, the present study identifies boundaries in the NLTVA parameter space delimiting safe, 
unsafe and unacceptable operations. The sensitivity of these boundaries to uncertainty in the NLTVA parameters is also investigated.

\vspace{1cm}

\end{abstract}
\vspace{-0.5cm}
\rule{0.85\linewidth}{.3pt}

\vspace{0.5cm} Corresponding author: \\ Thibaut Detroux\\
Space Structures and Systems Laboratory\\
Department of Aerospace and Mechanical Engineering\\
University of Liège
\\ 1 Chemin des Chevreuils (B52/3), B-4000 Liège, Belgium. \\
Email: tdetroux@ulg.ac.be
\vspace{2cm}\\}
\end{center}

\normalsize

\newpage

\section{Introduction}

A recent trend in the technical literature is to exploit nonlinear dynamical phenomena instead of avoiding them, as is the common practice. For instance, reference \cite{Daraio} demonstrates a new mechanism for tunable rectification that uses bifurcations and chaos. In \cite{Antonio}, a new strategy for engineering low-frequency noise oscillators is developed through the coupling of modes in internal resonance conditions. A cascade of parametric resonances is proposed by Strachan et al. as a basis for the development of passive frequency dividers \cite{Shaw}.

Nonlinearity is also more and more utilized for vibration absorption \cite{Poovarodom,Alexander,Machado,Lacarbonara} and energy harvesting \cite{Barton,Quinn,Worden,Inman}. For instance, a nonlinear energy sink (NES), i.e.,
an absorber with essential nonlinearity \cite{Book}, can extract energy from virtually any mode of a host structure \cite{NLD}. 
The NES can also carry out targeted energy transfer, which is an irreversible channeling of vibrational energy from the host structure to the absorber \cite{SIAM}. 
This absorber was applied for various purposes including seismic mitigation \cite{Nucera}, aeroelastic instability suppression \cite{Hubbard,Lamarque}, 
acoustic mitigation \cite{Cochelin} and chatter suppression \cite{Gourc}. Another recently-developed absorber is the nonlinear tuned vibration absorber (NLTVA) 
\cite{NLTVA1}. A unique feature of this device is that it can enforce equal peaks in the frequency response of the coupled system for a large range of motion 
amplitudes thereby generalizing Den Hartog's equal peak method to nonlinear systems. The NLTVA is therefore particularly suitable for mitigating the vibrations 
of a nonlinear resonance of a mechanical system. It was also found to be effective for the suppression of limit cycle oscillations \cite{ENOC}.

These contributions demonstrate that the purposeful introduction of nonlinearity can enhance system performance. However, nonlinearity can also give rise to complicated dynamical phenomena, which linear systems cannot. If quasiperiodic regimes of motion can be favorable for
vibration absorption with essential nonlinearity \cite{Staros}, they were found to be detrimental for a nonlinear absorber possessing both linear and
nonlinear springs \cite{Shaw4}. This highlights that no general conclusion can be drawn regarding the influence of quasiperiodic attractors.
Detached resonance curves (DRCs), also termed isolas, are generated by the multivaluedness of nonlinear responses and may limit the practical applicability of nonlinear absorbers \cite{staros2,Gourc2}. An important difficulty with DRCs is that they can easily be missed, because they are detached from the main resonance branch \cite{Alexander,Singh}. Finally, we note that DRCs were found in other applications involving nonlinearities, such as shimmying wheels \cite{Stepan} and structures with cyclic symmetry \cite{Grolet}, showing the generic character of DRCs.

In view of the potentially adverse effects of the aforementioned nonlinear attractors, the main objective of the present paper is to identify boundaries in the NLTVA parameter space delimiting safe, unsafe and unacceptable operations. The sensitivity of these boundaries to uncertainty in the NLTVA parameters is also investigated. To this end, rigorous nonlinear analysis methods, i.e., numerical continuation of periodic solutions, bifurcation detection and tracking, and global analysis, are utilized. Although these methods are well-established, their combination in a single study has not often been reported in the vibration mitigation literature.

The paper is organized as follows. Section 2 briefly reviews the salient features of the NLTVA. Specifically, this section demonstrates that equal peaks in the frequency response of the coupled system can be maintained in nonlinear regimes of motion. Section 3 reveals that systems featuring a NLTVA can exhibit DRCs and quasiperiodic regimes of motion. Based on the existence and location of these attractors, regions of safe, unsafe and unacceptable NLTVA operations are defined. Section 4 studies the sensitivity of attenuation performance and of the three regions of NLTVA operation to variations of the different absorber parameters. The conclusions of the present study are summarized in Section 5.

\section{Performance of the nonlinear tuned vibration absorber}

\begin{figure}[t]
\unitlength 1.3cm
\begin{picture}(10,2.5)(-7,1.5983)%
\Thicklines
\path(1.1016,3.8048)(2.7032,3.8048)(2.7032,2.2032)(1.1016,2.2032)(1.1016,3.8048)
\path(-2.336,3.8048)(-0.7344,3.8048)(-0.7344,2.2032)(-2.336,2.2032)(-2.336,3.8048)

\allinethickness{0.4pt}%
\path(-4.172,3.0294)(-3.8966,3.0294)(-3.8048,3.213)(-3.6212,2.8458)(-3.4376,3.213)(-3.254,2.8458)(-3.0704,3.213)
(-2.8868,2.8458)(-2.7032,3.213)(-2.6114,3.0294)(-2.336,3.0294)

\path(-4.172,2.4786)(-3.8966,2.4786)(-3.8048,2.6622)(-3.6212,2.295)(-3.4376,2.6622)(-3.254,2.295)(-3.0704,2.6622)
(-2.8868,2.295)(-2.7032,2.6622)(-2.6114,2.4786)(-2.336,2.4786)
\path(-3.2735,2.6805)(-3.1622,2.754)(-3.1542,2.6209)
\path(-3.4376,2.2032)(-3.1622,2.754)

\path(-0.7344,2.75)(-0.45,2.75)(-0.45,2.5)(0.82,2.5)(0.82,3.0)(-0.45,3.0)(-0.45,2.75)
\path(0.82,2.75)(1.1016,2.75)
\put(0.2,2.75){\makebox(0,0){$?$}}

\path(-4.172,3.57)(-3.45,3.57)\path(-3.15,3.57)(-2.336,3.57)
\path(-3,3.35)(-3.45,3.35)(-3.45,3.79)(-3,3.79)
\path(-3.15,3.42)(-3.15,3.72)

\path(-0.736,3.57)(-0.014,3.57)\path(0.286,3.57)(1.10,3.57)
\path(0.436,3.35)(-0.014,3.35)(-0.014,3.79)(0.436,3.79)
\path(0.286,3.42)(0.286,3.72)

\thicklines
\path(-4.172,4.01)(-4.172,1.99)
\allinethickness{0.4pt}
\path(-4.549,4.008)(-4.631,3.9268)\path(-4.3369,4.008)(-4.631,3.7147)\path(-4.172,3.9615)(-4.631,3.5025)
\path(-4.172,3.7494)(-4.631,3.2904)\path(-4.172,3.5373)(-4.631,3.0783)\path(-4.172,3.3251)(-4.631,2.8661)
\path(-4.172,3.113)(-4.631,2.654)\path(-4.172,2.9009)(-4.631,2.4419)\path(-4.172,2.6887)(-4.631,2.2297)
\path(-4.172,2.4766)(-4.631,2.0176)\path(-4.172,2.2645)(-4.4473,1.9892)\path(-4.172,2.0523)(-4.2351,1.9892)

\path(-1.53,2.1801)(-1.53,1.6157)
\path(-1.28,1.88)(-1.18,1.82)(-1.28,1.76)
\path(-1.53,1.82)(-1.18,1.82)
\put(-0.9,1.7921){\makebox(0,0){$x_1$}}

\path(2.15,1.894)(2.25,1.8273)(2.15,1.7607)
\path(1.9,1.8273)(2.25,1.8273)
\path(1.9,2.2154)(1.9,1.6509)
\put(2.53,1.79){\makebox(0,0){$x_2$}}

\put(-1.5,3){\makebox(0,0){$m_1$}}
\put(1.9,3){\makebox(0,0){$m_2$}}

\put(-3.2,4.){\makebox(0,0){$c_1$}}
\put(-2.6,3.35){\makebox(0,0){$k_1$}}
\put(-2.6,2.1){\makebox(0,0){$k_{nl1}$}}
\put(0.2,4.){\makebox(0,0){$c_2$}}
\put(0.2,2.2){\makebox(0,0){$g(\bullet)$}}

\end{picture}%
\caption{Schematic representation of an NLTVA attached to a Duffing oscillator.}\label{2dofsystem}
\end{figure}
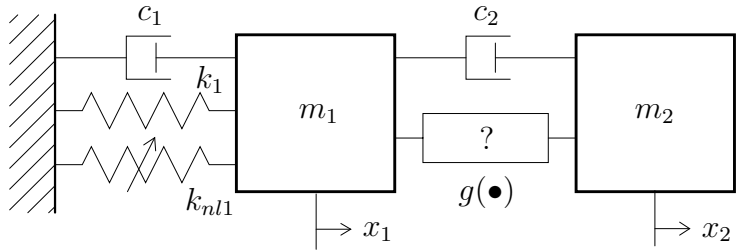

The NLTVA targets the mitigation of a nonlinear resonance in an as large as possible range of forcing amplitudes. An unconventional feature
of this absorber is that the mathematical form of its nonlinear restoring force is not imposed a priori, as it is the case for most
existing nonlinear absorbers. Instead, we fully exploit the additional design parameter offered by nonlinear devices, and, hence, we synthesize
the absorber's load-deflection curve according to the nonlinear restoring force of the primary structure.

The dynamics of a Duffing oscillator with an attached NLTVA, as depicted in Figure \ref{2dofsystem}, is considered throughout this study:
\begin{eqnarray}\label{eomdim}
   \nonumber
   m_1\ddot{x}_1+c_1\dot{x}_1+k_1x_1+k_{nl1}x_1^3+c_2(\dot{x}_1-\dot{x}_2)+g(x_1-x_2)& = & F\cos \omega t\\
   m_2\ddot{x}_2+c_2(\dot{x}_2-\dot{x}_1)-g(x_1-x_2) & = & 0
\end{eqnarray}
where $x_1(t)$ and $x_2(t)$ are the displacements of the harmonically-forced primary system and of the NLTVA, respectively. 
The NLTVA is assumed to have a generic smooth restoring force $g\left(x_1-x_2\right)$ with $g(0)=0$. 
In order to avoid important sensitivity of absorber performance to forcing amplitude, it was shown in reference \cite{NLTVA1} that the 
function $g(x_1-x_2)$ should be chosen such that the NLTVA is a `mirror' of the primary system. More precisely, 
besides a linear spring, the NLTVA should possess a nonlinear spring of the same mathematical form as that of the nonlinear spring of the primary system. 
For instance, if the nonlinearity in the primary system is quadratic or cubic, the NLTVA should possess a quadratic or a cubic spring, respectively. 
To mitigate the vibrations of the Duffing oscillator, a NLTVA with linear and cubic stiffnesses is therefore considered, 
i.e., $g(x_1-x_2)=k_2(x_1-x_2)+k_{nl2}(x_1-x_2)^3$, and the governing equations of motion becomes
\begin{eqnarray}
   \nonumber
   m_1\ddot{x}_1+c_1\dot{x}_1+k_1x_1+k_{nl1}x_1^3+c_2(\dot{x}_1-\dot{x}_2)+k_2(x_1-x_2)+k_{nl2}(x_1-x_2)^3& = & F\cos \omega t\\
   m_2\ddot{x}_2+c_2(\dot{x}_2-\dot{x}_1)+k_2(x_2-x_1)+k_{nl2}(x_2-x_1)^3 & = & 0\label{eomdim2}
\end{eqnarray}

In view of the effectiveness of the equal-peak method \cite{DenHartog,Brock} for the design of linear tuned vibration absorbers (LTVA) attached to linear host structures, an attempt to generalize this tuning rule to nonlinear absorbers attached to nonlinear host structures was made in reference \cite{NLTVA1}. The first step was to impose equal peaks in the receptance function of the underlying linear system using the formulas proposed by Asami et al. \cite{Asami}:
\begin{eqnarray}\label{DHrule_dim}
   \nonumber k_2^{opt}&=&\displaystyle\frac{8\epsilon k_1\left[16+23\epsilon+9\epsilon^2+2(2+\epsilon)\sqrt{4+3\epsilon} \right]}
    {3(1+\epsilon)^2(64+80\epsilon+27\epsilon^2)}\\
   c_2^{opt}&=&\displaystyle\sqrt{\frac{k_2m_2(8+9\epsilon-4\sqrt{4+3\epsilon})}{4(1+\epsilon)}}
  \end{eqnarray}
where $\epsilon = m_2/m_1$ is the mass ratio, chosen according to practical constraints. We note that these formulas are exact for an undamped primary system, unlike those proposed previously by Den Hartog \cite{DenHartog} and Brock \cite{Brock}.
The second step was to determine the nonlinear coefficient $k_{nl2}$ that can maintain two resonance peaks of equal amplitude in nonlinear regimes of motion. 
A very interesting result of \cite{NLTVA1} is that the nonlinear coefficient that realizes equal peaks for various forcing amplitudes is almost constant and can be accurately calculated using the analytical expression:
\begin{eqnarray}\label{DHrule_dim2}
    k_{nl2}^{opt}&=&\displaystyle\frac{2\epsilon^2k_{nl1}}{(1+4\epsilon)}
\end{eqnarray}
Overall, Equations (\ref{DHrule_dim}) and (\ref{DHrule_dim2}) represent a new tuning rule for nonlinear vibration absorbers that may be viewed as a nonlinear generalization of Den Hartog's equal-peak method.

\begin{figure}[t]
\setlength{\unitlength}{1cm}
\begin{picture}(8,8.5)(0,0)
\put(3,0.6){\includegraphics[width=10truecm]{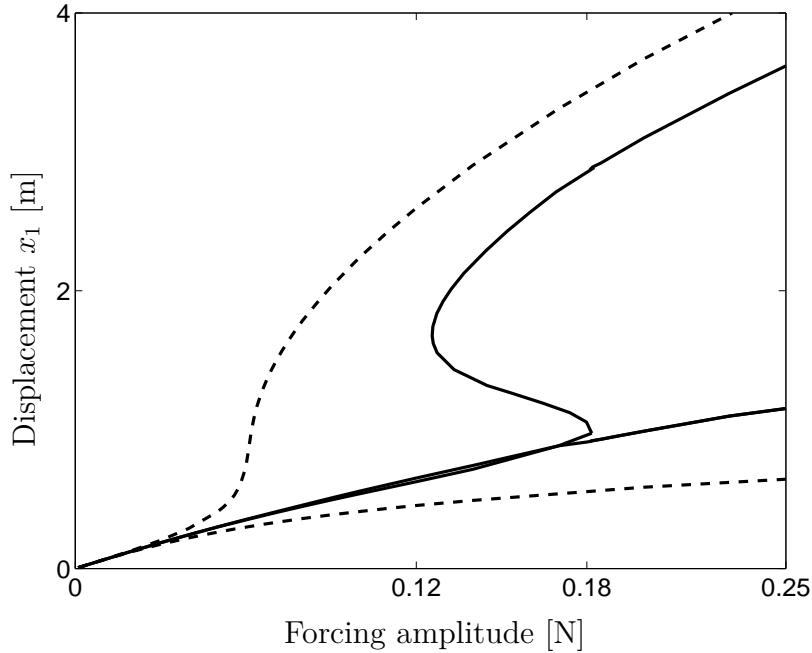}}
\put(6.0,0){Forcing amplitude [N]} \put(2.4,2.6){\rotatebox{90}{Displacement $x_1$ [m]}}
\end{picture}\caption{Performance of the LVTA/NLTVA for increasing forcing amplitudes. The dashed and solid lines depict the amplitude
of the resonances peaks of the Duffing oscillator with an attached LTVA and NLTVA, respectively.}\label{compa_NLTVA}
\end{figure}

A comparison of the performance of a LTVA and of a NLTVA attached to a Duffing oscillator is performed for the parameters listed in Table \ref{tab:param_sys}. Figure \ref{compa_NLTVA} represents the amplitude of the resonance peaks for increasing forcing amplitudes. The first observation is that the NLTVA performance (in terms of $H_\infty$ optimization) is always superior to that of
the LTVA. For the LTVA, the two resonance peaks start to have different amplitudes from $F = 0.03\,$N, showing the detuning of this absorber in nonlinear regimes of motion. Conversely, for the NLTVA, the amplitudes of the two resonance peaks remain almost identical until $F = 0.18\,$N, providing the numerical evidence of the effectiveness of the design proposed in Equations (\ref{DHrule_dim}) and (\ref{DHrule_dim2}). However, between $F = 0.12\,$N and $0.18\,$N, the two main resonance peaks co-exist with two additional resonance peaks that will be shown to correspond to a DRC.

\begin{table}[b]
\begin{center}
\begin{tabular}{c||c|c|c}
 & Primary system & LTVA & NLTVA\\
\hline
\hline
Mass [kg] & $m_{1} = 1$ & $m_{2} = 0.05$ & $m_{2} = 0.05$\\

Linear stiffness [N/m] & $k_{1} = 1$ & $k_{2} = 0.0454$ & $k_{2} = 0.0454$\\

Linear damping [Ns/m] & $c_{1} = 0.002$ & $c_{2} = 0.0128$ & $c_{2} = 0.0128$\\

Nonlinear stiffness [N/m$^3$] & $k_{nl1} = 1$ & --- & $k_{nl2} = 0.0042$\\
\end{tabular}
\caption{Parameters of the Duffing oscillator and attached LTVA and NLTVA ($\epsilon=0.05$).}
\label{tab:param_sys}
\end{center}
\end{table}

\section{Robustness of the nonlinear tuned vibration absorber}\label{intro_dyn}

Considering the multivaluedness of the NLTVA response in Figure \ref{compa_NLTVA}, this section aims at uncovering the dynamical attractors that a Duffing oscillator 
featuring a NLTVA may exhibit. The methodology in this study combines numerical continuation of periodic solutions, 
bifurcation detection and tracking, and global analysis.

\subsection{Numerical continuation of periodic solutions and detection of bifurcations}

The frequency response of system (\ref{eomdim2}) for the parameters listed in Table \ref{tab:param_sys} was computed using the algorithm proposed in reference \cite{Orlando}. Codimension-1 numerical continuation was carried out using the multi-harmonic balance method, which approximates periodic solutions using Fourier series. Stability analysis was achieved using Hill's method, and fold and Neimark-Sacker (NS) bifurcations were detected using appropriate test functions. Convergence of the results was obtained when the first 5 harmonics were retained. 

\begin{figure}[p]
\setlength{\unitlength}{1cm}
\begin{picture}(8,18)(0,0)
\put(0.75,13.4){\includegraphics[width=7.0truecm]{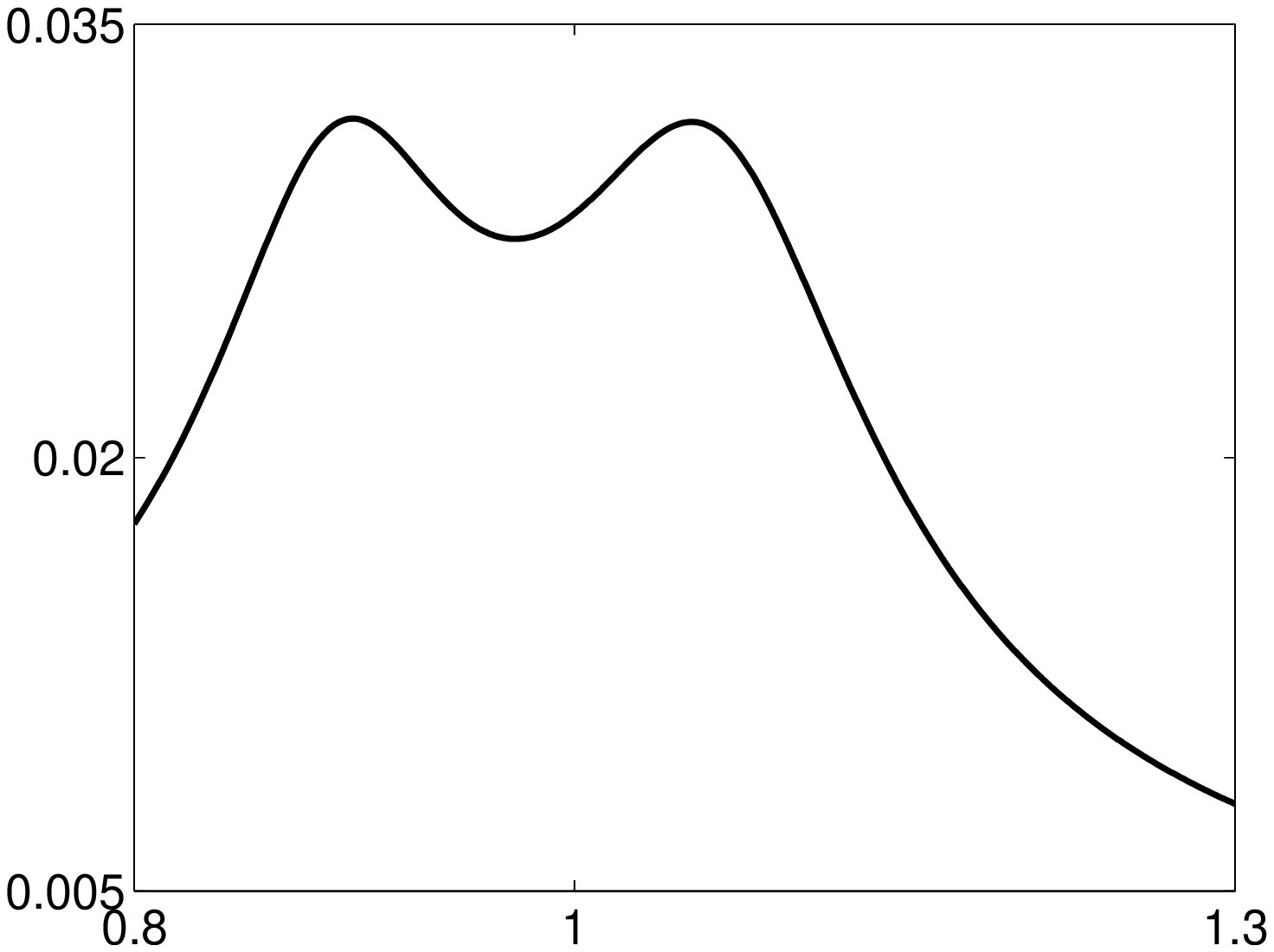}}
\put(8.45,13.4){\includegraphics[width=7.0truecm]{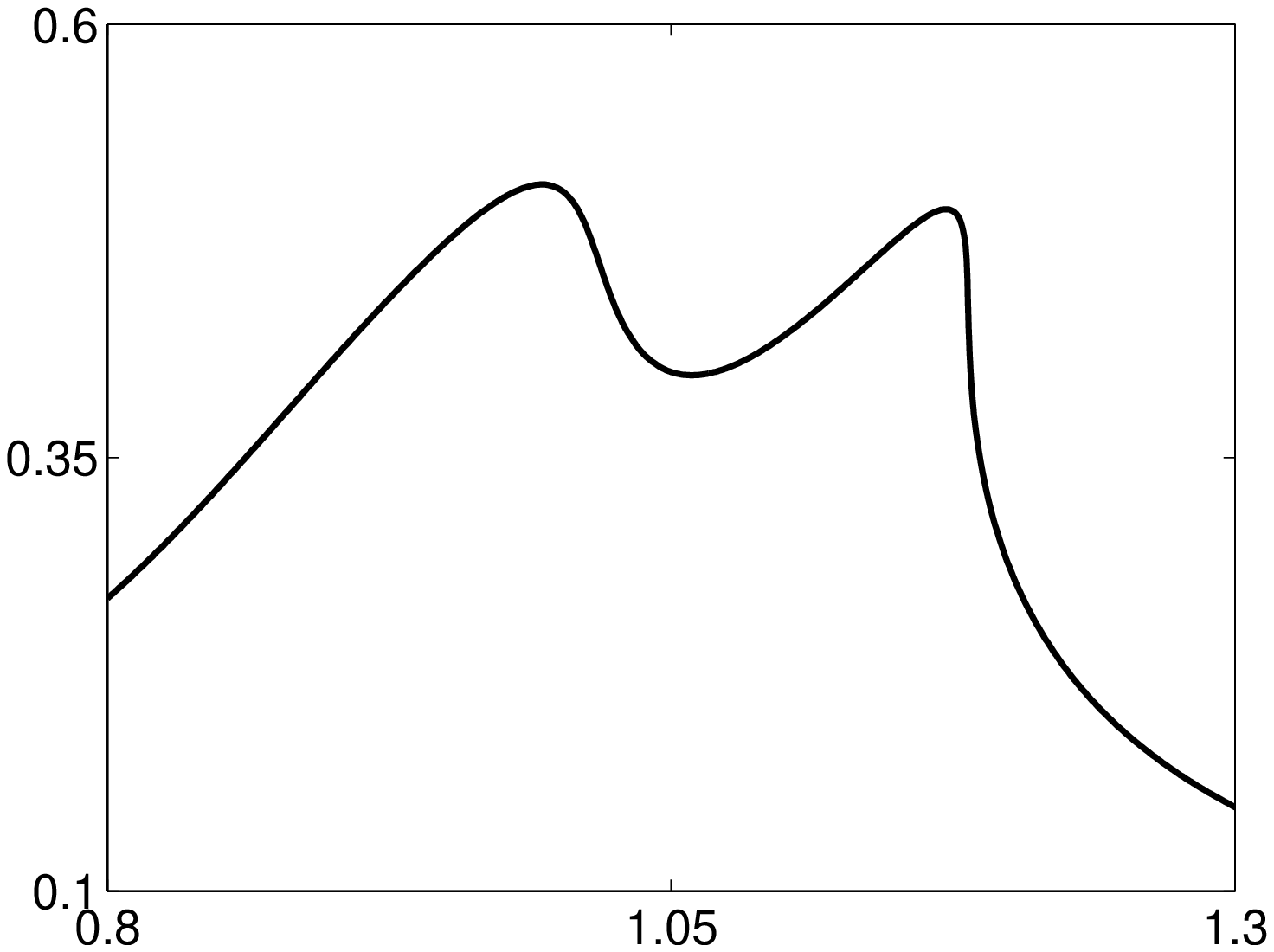}}
\put(0.75,7){\includegraphics[width=7.0truecm]{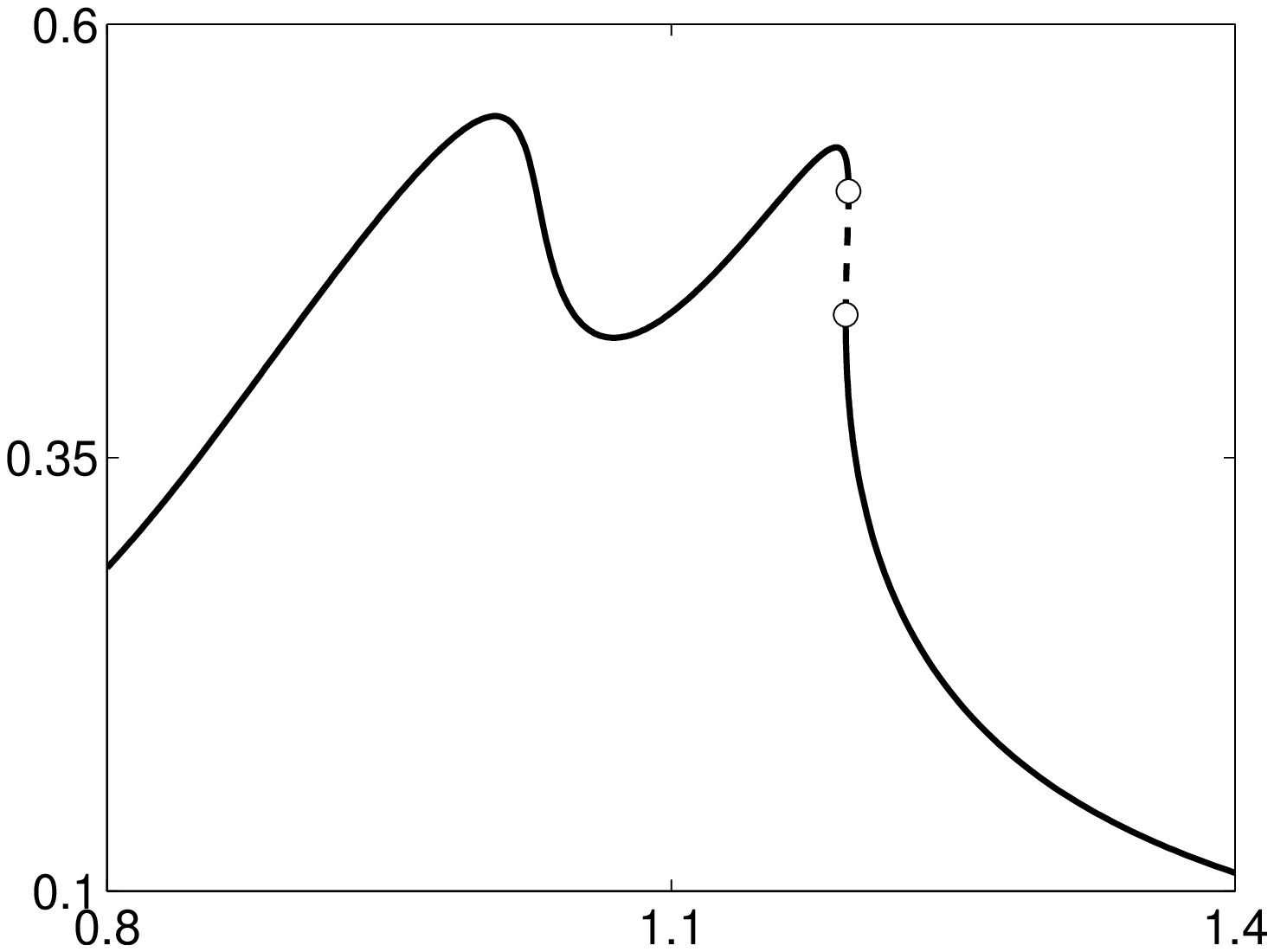}}
\put(8.45,7){\includegraphics[width=7.0truecm]{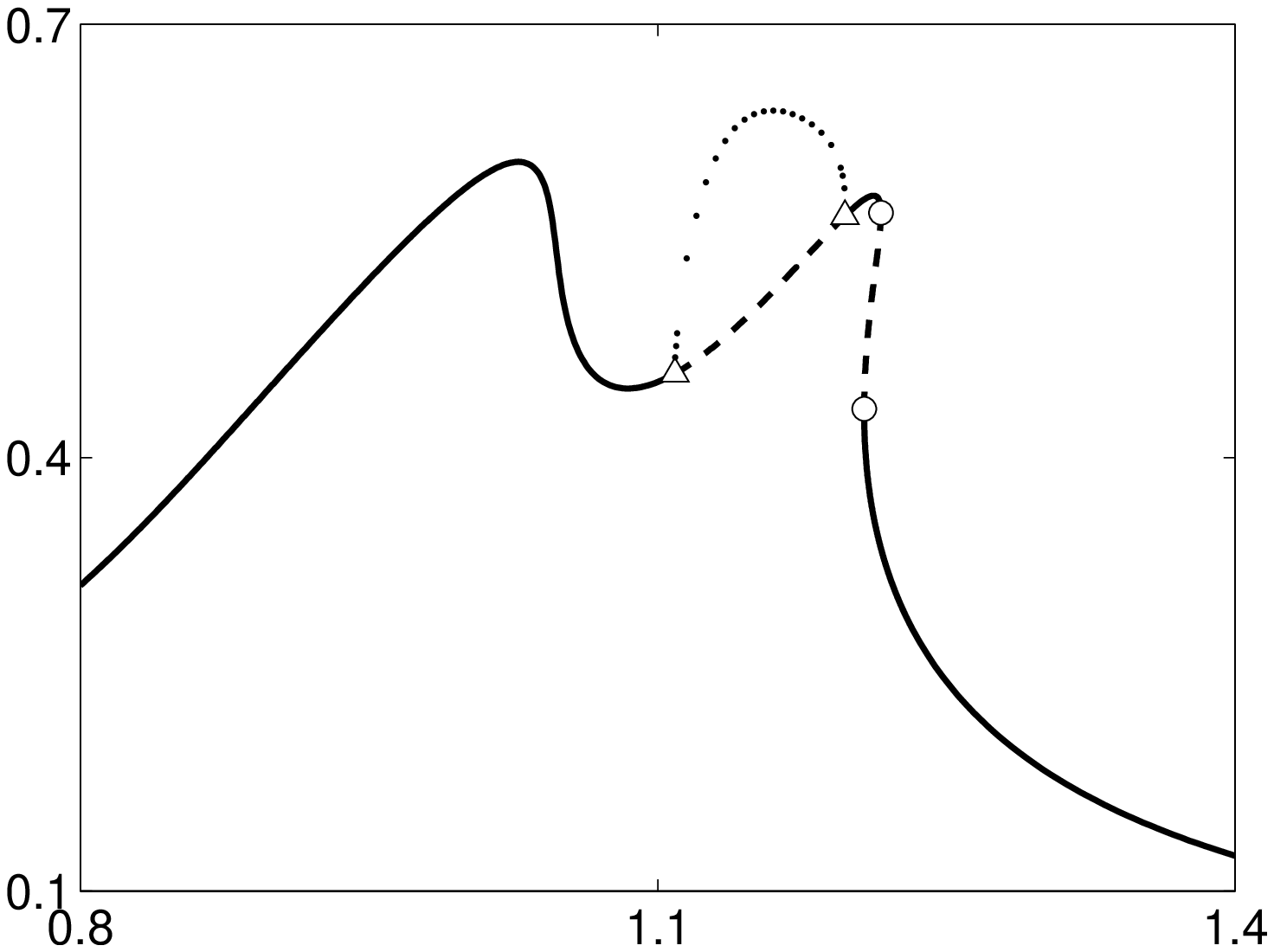}}
\put(0.75,0.6){\includegraphics[width=7.0truecm]{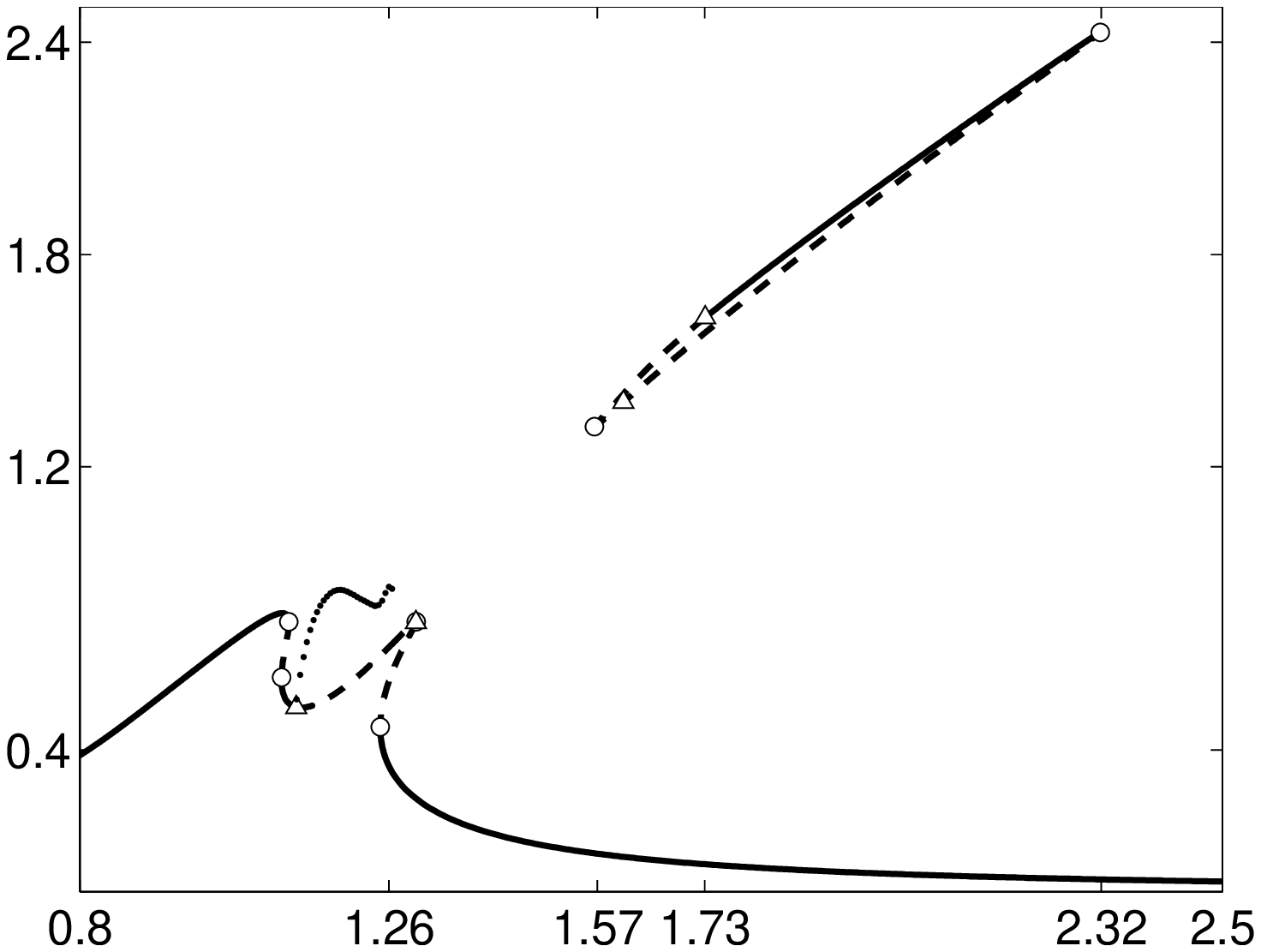}}
\put(8.45,0.6){\includegraphics[width=7.0truecm]{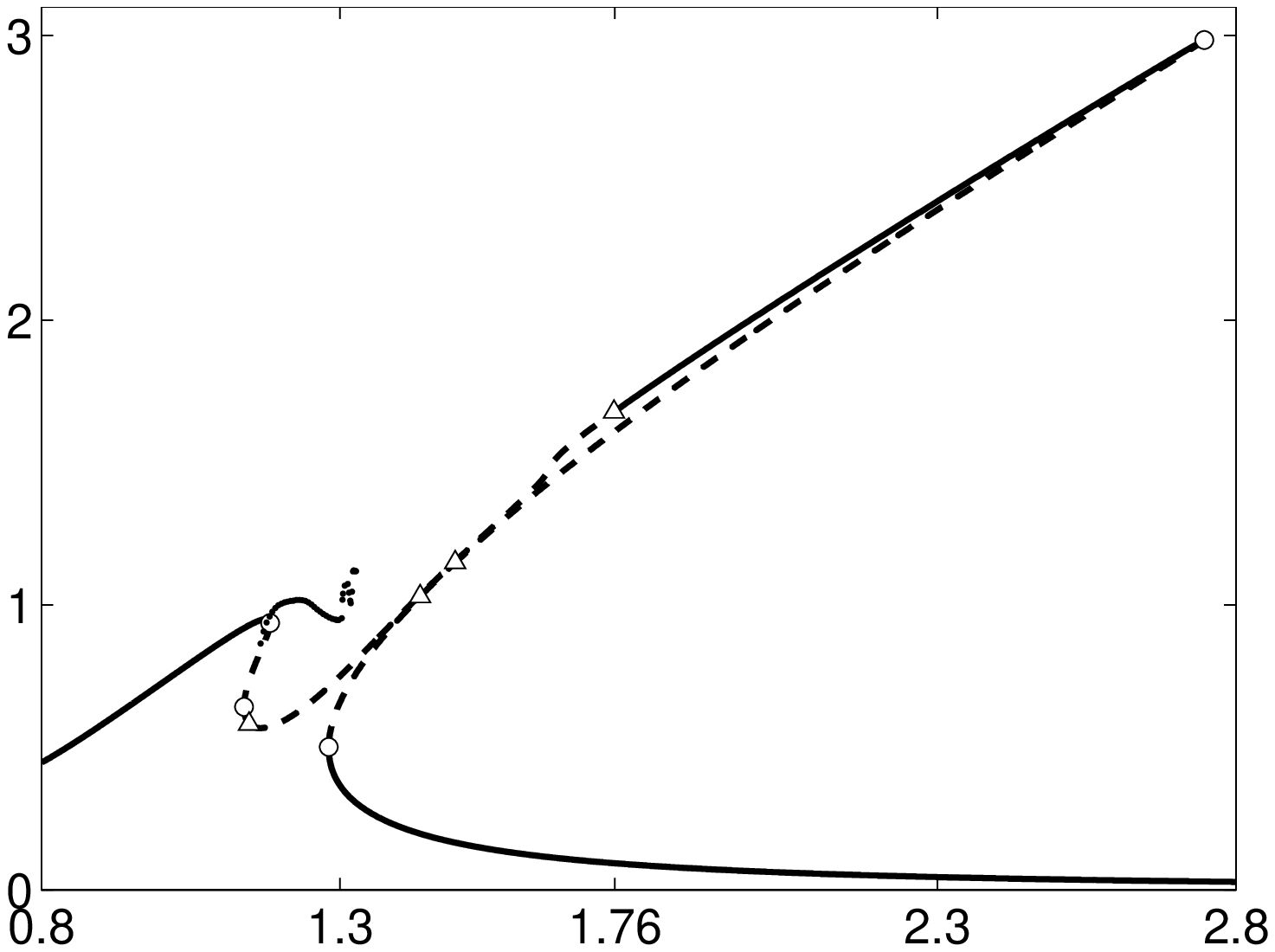}}
\put(6.5,0){Frequency [rad/s]} \put(-0.1,7.9){\rotatebox{90}{Displacement $x_1$ [m]}}
 \put(1.5,18.0){(a)}
 \put(9.1,18.0){(b)}
 \put(1.5,11.6){(c)}
 \put(9.1,11.6){(d)}
 \put(1.5,5.2){(e)}
 \put(9.1,5.2){(f)}
\end{picture}\caption{Frequency response of the Duffing oscillator with an attached NLTVA. (a) $F = 0.005\,$N; (b) $F=0.09\,$N; (c) $F = 0.098\,$N; (d) $F=0.11\,$N; (e) $F=0.15\,$N; (f) $F=0.19\,$N. The solid and dashed lines represent stable and unstable solutions, respectively. Fold and Neimark-Sacker bifurcations are depicted with circle and triangle markers, respectively. The dotted line represents stable quasiperiodic oscillations.}\label{Freq_resp_regimes}
\end{figure}

For $F = 0.005\,$N in Figure \ref{Freq_resp_regimes}(a), neither of the nonlinearities are activated, and the classical linear result is retrieved \cite{DenHartog}. For $F = 0.09\,$N in Figure  \ref{Freq_resp_regimes}(b), the resonance peaks bend forward as a result of the hardening nature of the cubic springs, and the resonance frequencies increase. Notwithstanding this nonlinear behavior, resonance peaks of equal amplitude are obtained thanks to the NLTVA.

Slightly increasing forcing amplitude triggers the appearance of two different bifurcations. For $F = 0.098\,$N, a pair of fold bifurcations modifies the stability along the frequency response in Figure \ref{Freq_resp_regimes}(c). For $F = 0.11\,$N in Figure \ref{Freq_resp_regimes}(d), a pair of NS bifurcations changes stability as well, but it also generates a stable branch of quasiperiodic solutions that was computed using direct time integrations. Since approximately equal peaks are maintained in Figures \ref{Freq_resp_regimes}(c-d) and since quasiperiodic oscillations have amplitudes comparable to those of the resonance peaks, the NLTVA can still be considered as effective.

For $F = 0.15\,$N, a pair of fold bifurcations creates a DRC between 1.57 and 2.32 rad/s in Figure \ref{Freq_resp_regimes}(e). This DRC is associated with large amplitudes of motion, but it remains far enough from the desired operating frequency range of the NLTVA. In addition, the left portion of the DRC between 1.57 and 1.73 rad/s is unstable, and, hence, not physically realizable. The DRC also possesses a pair of NS bifurcations, but no stable branch of quasiperiodic oscillations could be found. 

For larger forcing amplitudes, the DRC expands and eventually merges with the second resonance peak for $F = 0.19\,$N, as depicted in Figure \ref{Freq_resp_regimes}(f). This merging eliminates the fold bifurcation characterizing the second resonance peak and the fold bifurcation on the left of the DRC, and causes a very substantial increase in the amplitude of the second resonance.

\subsection{Bifurcation tracking}

Since the creation of the quasiperiodic solutions and of the DRC occurs through NS and fold bifurcations, respectively, these bifurcations were tracked in the three-dimensional space ($x_1,\omega,F$) using the multi-harmonic balance method \cite{Orlando}. To facilitate the interpretation of the results, the bifurcation loci were projected onto the two-dimensional plane ($x_1,F$).

\begin{figure}[t]
\setlength{\unitlength}{1cm}
\begin{picture}(8,8.5)(0,0)
\put(3,0.6){\includegraphics[width=10truecm]{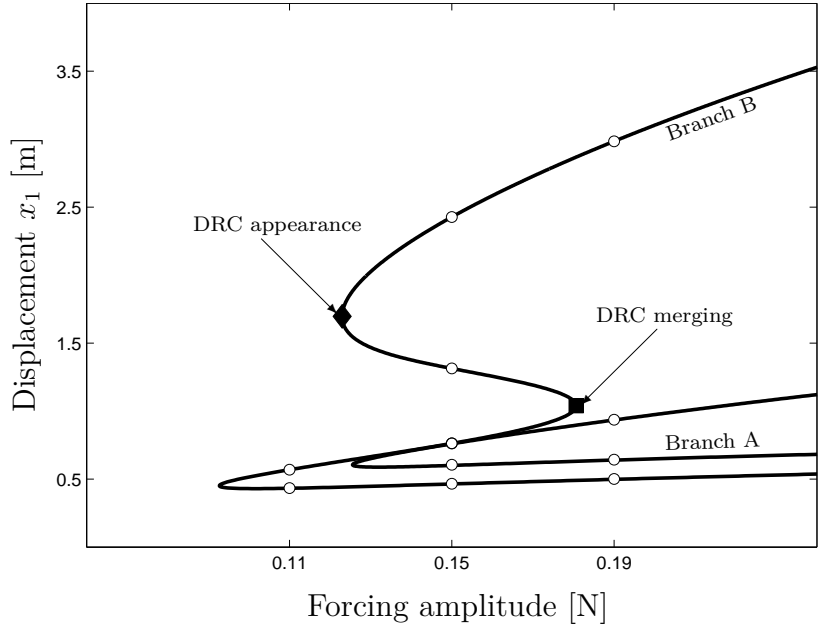}}
\put(6.3,0){Forcing amplitude [N]} \put(2.4,2.7){\rotatebox{90}{Displacement $x_1$ [m]}}
\put(11,2.2){\rotatebox{2}{\scriptsize{Branch A}}}
\put(11,6.3){\rotatebox{20}{\scriptsize{Branch B}}}
\put(4.8,5.1){\scriptsize{DRC appearance}}\put(5.7,5.0){\vector(1,-1){1}}
\put(10.1,3.9){\scriptsize{DRC merging}}\put(10.9,3.8){\vector(-1,-1){1}}
\end{picture}\caption{Projection of the branches of fold bifurcations onto the ($x_1,F$) plane. The circle markers represents the fold bifurcations in Figures \ref{Freq_resp_regimes}(d-f). The diamond and square markers indicate the appearance and merging of the DRCs, respectively.}\label{LP_tracking_analysis}
\end{figure}

Figure \ref{LP_tracking_analysis} represents the loci of the fold bifurcations of Figure \ref{Freq_resp_regimes}. Branch A is related to the bifurcations in the neighbourhood of the first resonance peak, whereas branch B corresponds to the bifurcations in the vicinity of the second resonance peak. Because the DRC merges with this latter peak, branch B indicates the creation and elimination of the DRC, represented with diamond and square markers, respectively. It can therefore be concluded that the DRC appears around $0.12\,$N and merges with the main branch around $0.18\,$N. It can also be observed that branches A and B overlap between 0.13 and 0.17 N, which is another manifestation of the proposed nonlinear equal-peak method. 

\begin{figure}[h]
\setlength{\unitlength}{1cm}
\begin{picture}(8,8.5)(0,0)
\put(3,0.6){\includegraphics[width=10truecm]{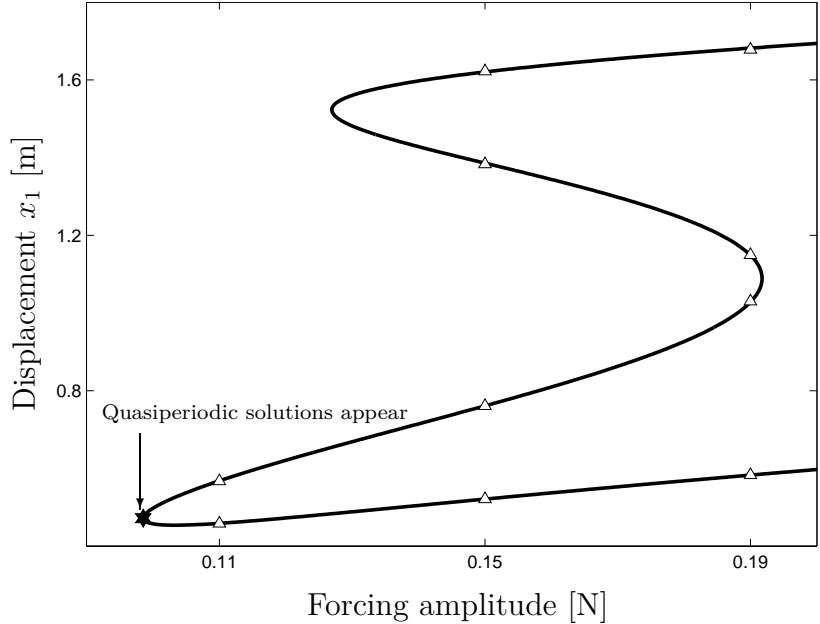}}
\put(3.6,2.6){\scriptsize{Quasiperiodic solutions appear}}\put(4.1,2.4){\vector(0,-1){1}}
\put(6.3,0){Forcing amplitude [N]} \put(2.4,2.7){\rotatebox{90}{Displacement $x_1$ [m]}}
\end{picture}\caption{Projection of the branch of Neimark-Sacker bifurcations onto the ($x_1,F$) plane. The triangle markers represent Neimark-Sacker bifurcations in Figures \ref{Freq_resp_regimes}(d-f). The star marker indicates the appearance of quasiperiodic solutions.}\label{NS_tracking_analysis}
\end{figure}

Figure \ref{NS_tracking_analysis} displays the locus of NS bifurcations and indicates that the stable NS branch between the two resonance peaks is created when $F = 0.095\,$N (depicted with a star marker). Interestingly, because the DRC possesses NS bifurcations in addition to the fold bifurcations, 
the existence of the DRC is also revealed by the upper turning point in Figure \ref{NS_tracking_analysis}.

\subsection{Global analysis of the adverse dynamics}\label{global_anal}
 
The previous two sections have highlighted the different dynamical attractors of the coupled system (\ref{eomdim2}) together with their locations in the ($F,\omega$) plane. Since the stable quasiperiodic solutions between the two resonances are associated with 
acceptable amplitudes, they are not considered as a great concern. However, we should stress that NS bifurcations 
can trigger torus breakdown and phase locking with large-amplitude attractors \cite{Afrai}. Despite detailed 
numerical simulations and co-dimension 2 bifurcation analysis in MatCont \cite{Matcont}, no evidence of such behaviors could be observed.

The large-amplitude DRC is more problematic. The likelihood of converging to the safe, low-amplitude periodic solution in the region where the DRC exists should therefore be determined. Given a forcing amplitude and frequency, direct time integrations for a large set of random initial states provided the basins of attraction, similarly to what was achieved in \cite{Dick} for a coupled linear oscillator and nonlinear absorber system. To limit the scope of the discussion, the NLTVA was considered at rest, as it would be the case, e.g., during an earthquake, but other configurations were also tested.

\begin{figure}[p]
\setlength{\unitlength}{1cm}
\begin{picture}(8,18)(0,0)
\put(0.75,13.4){\includegraphics[scale = 0.5]{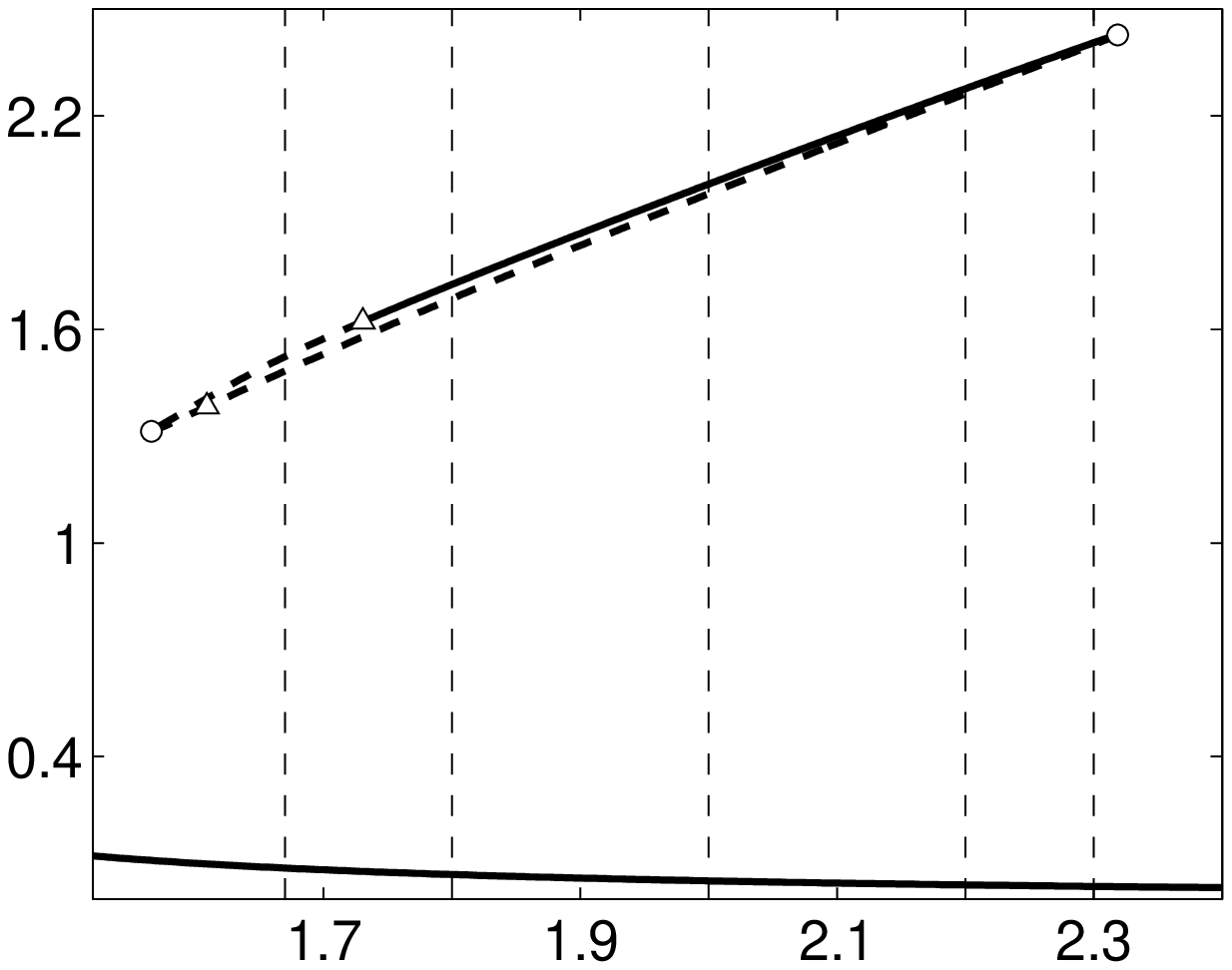}}
\put(8.65,13.4){\includegraphics[scale = 0.5]{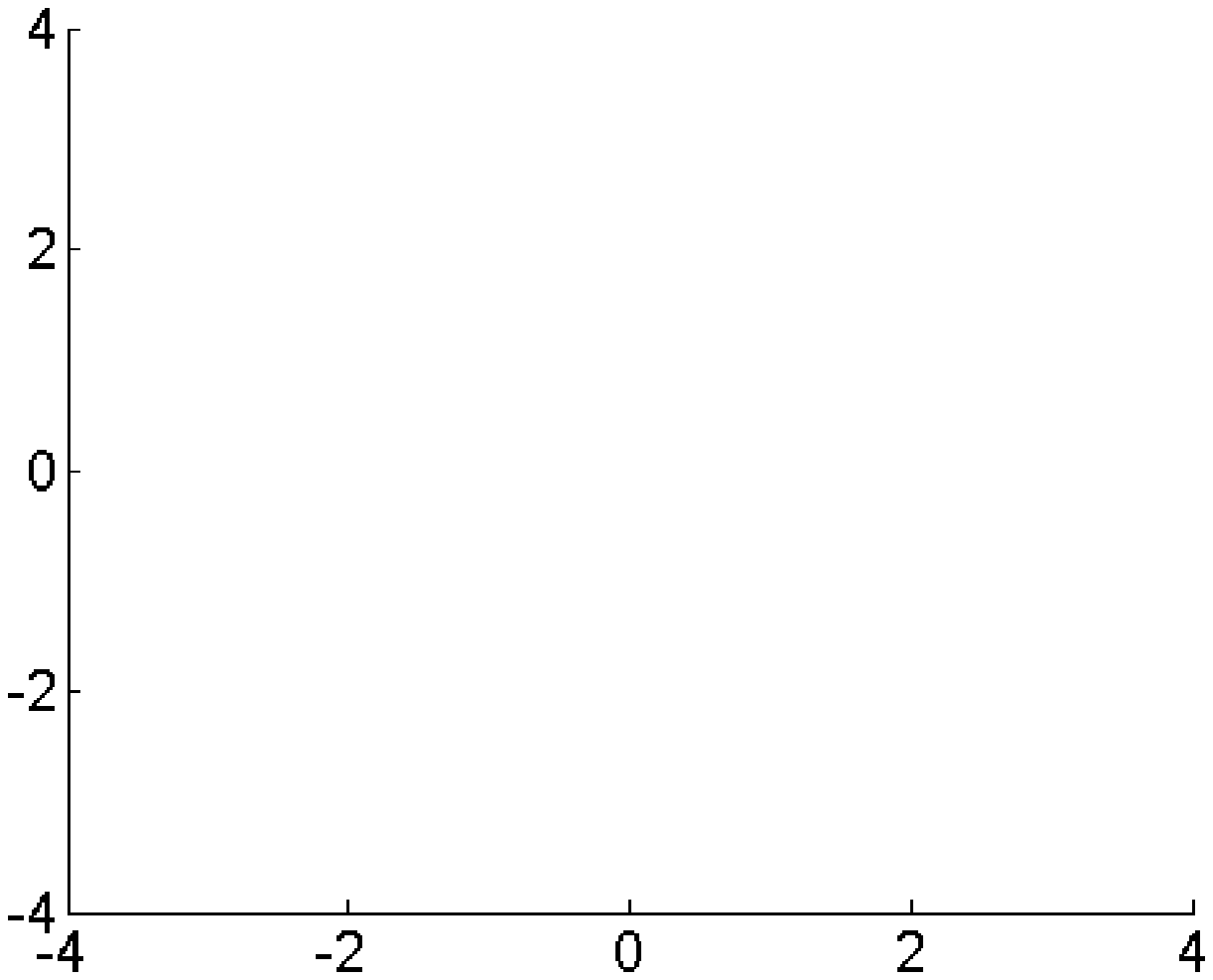}}
\put(0.75,7){\includegraphics[scale = 0.5]{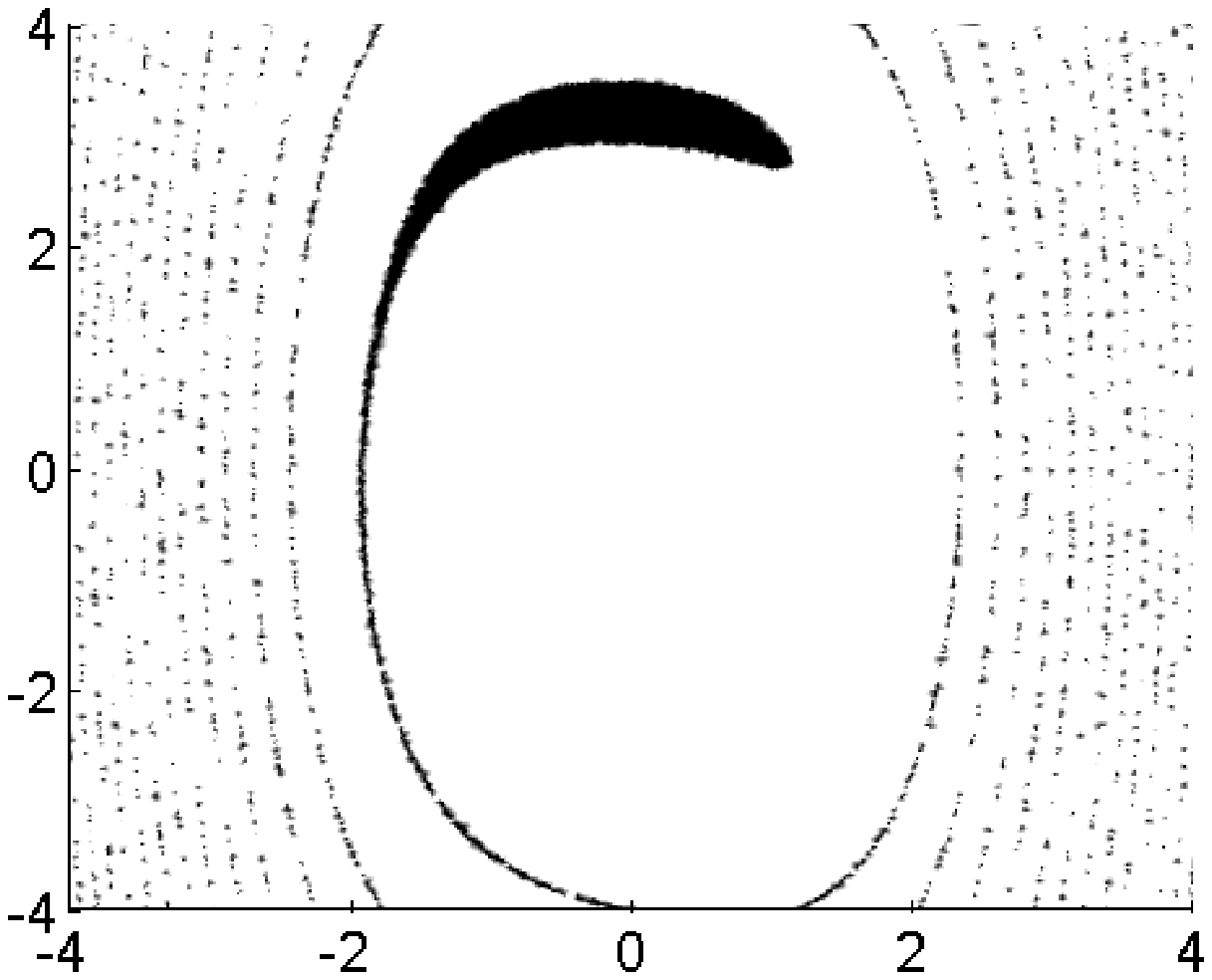}}
\put(8.65,7){\includegraphics[scale = 0.5]{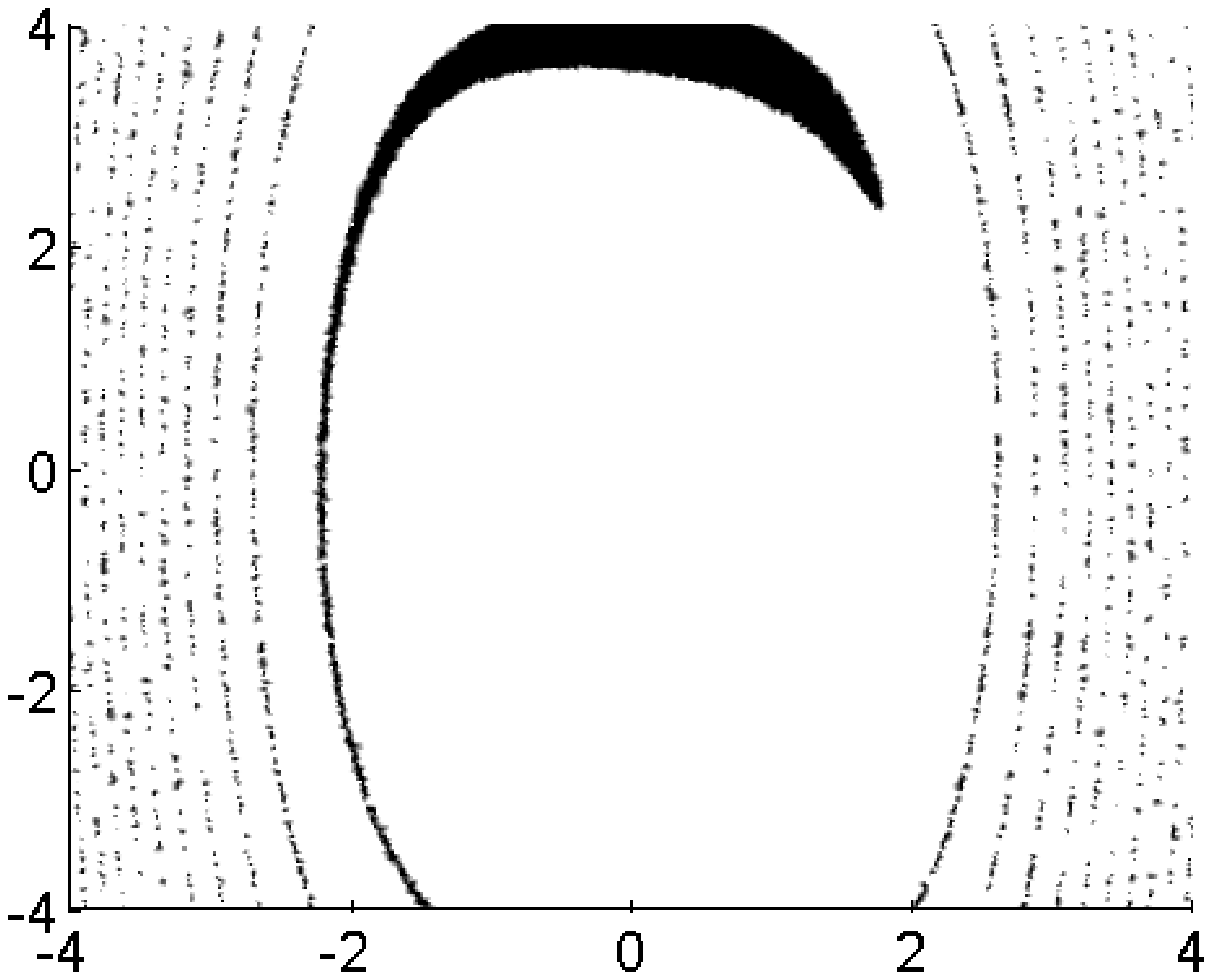}}
\put(0.75,0.6){\includegraphics[scale = 0.5]{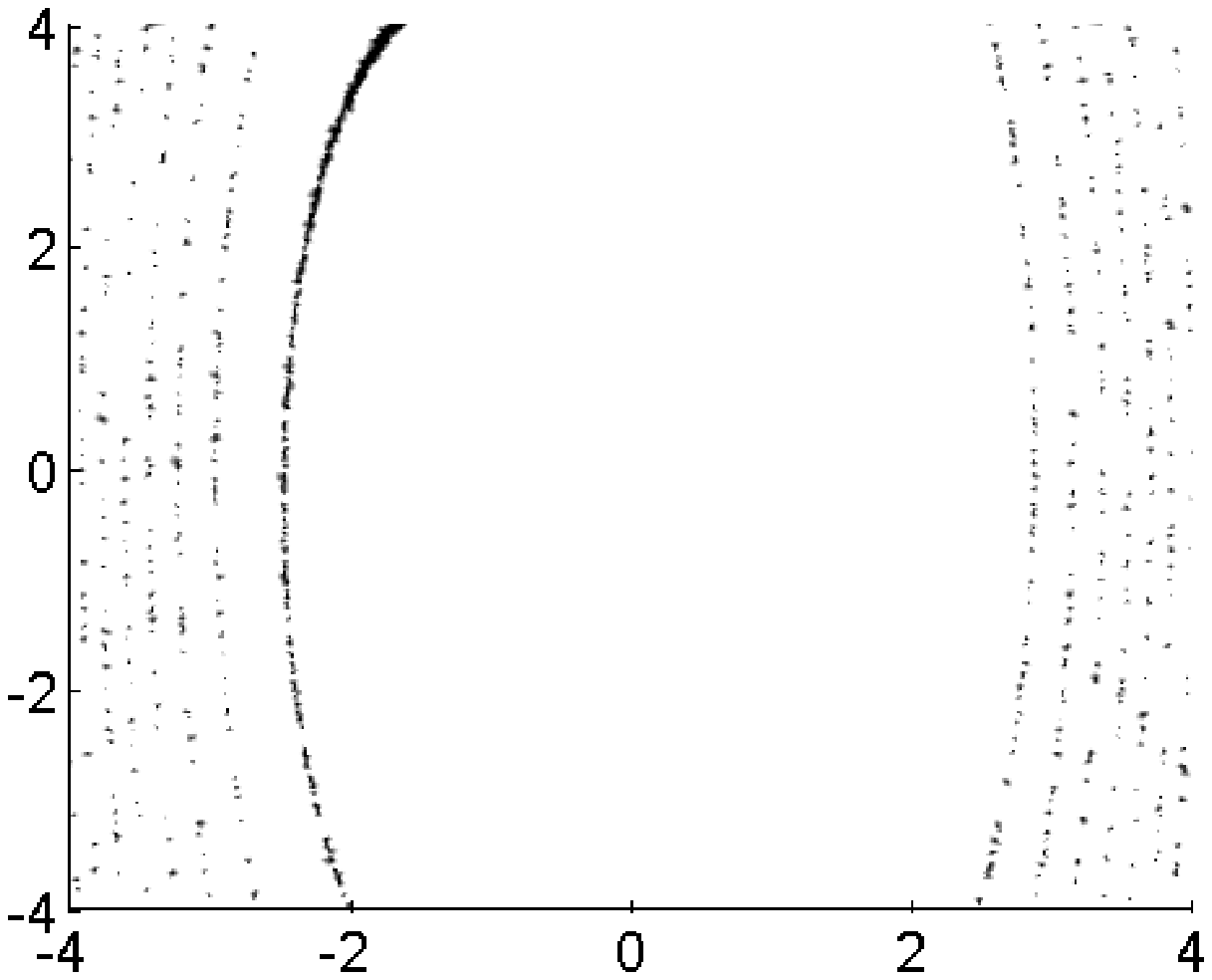}}
\put(8.65,0.6){\includegraphics[scale = 0.5]{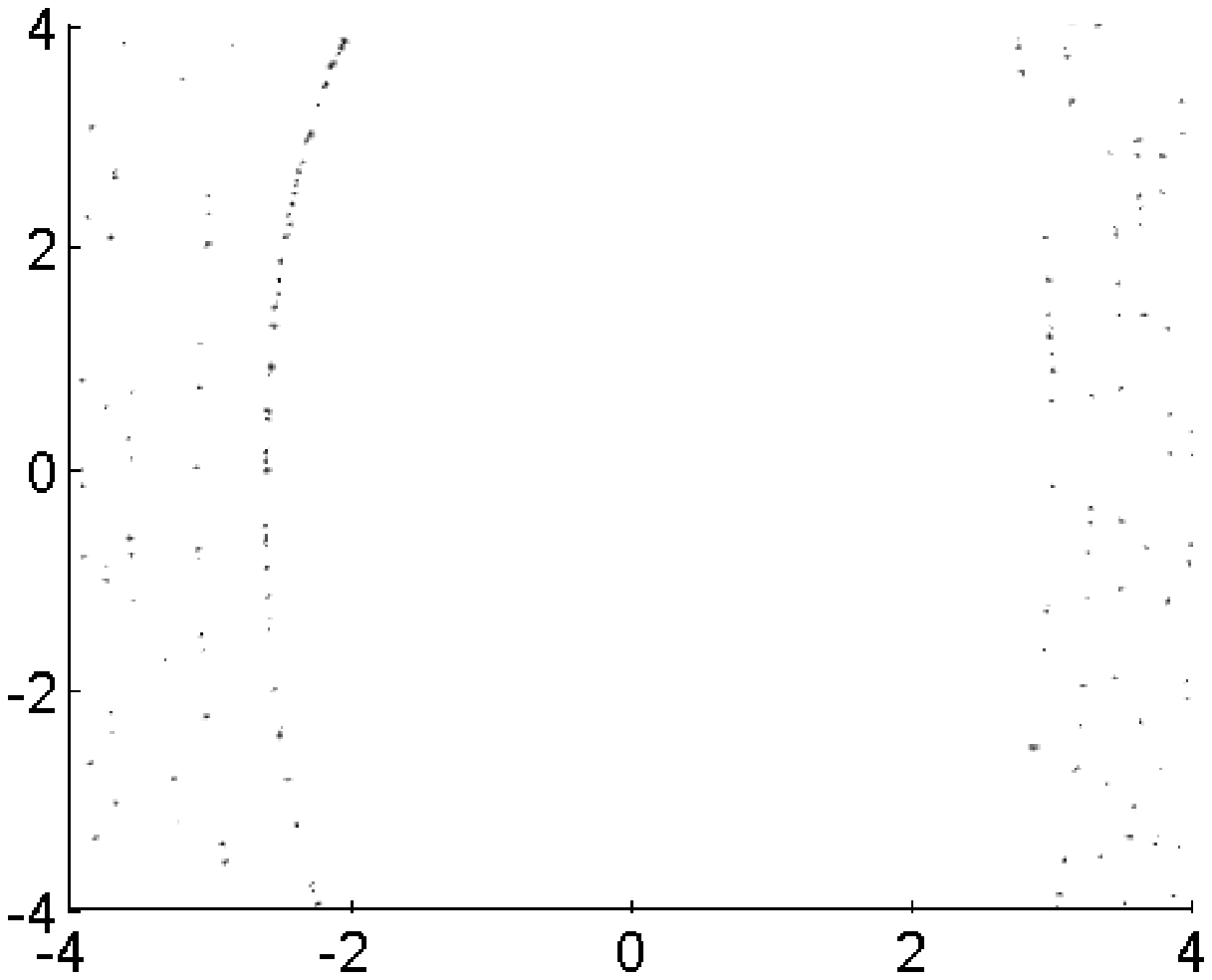}}
\put(3.0,13.2){Frequency [rad/s]}
\put(2.88,6.8){Initial displ. $x_1^0$ [m]}
\put(2.88,0.4){Initial displ. $x_1^0$ [m]}
\put(10.78,13.2){Initial displ. $x_1^0$ [m]}
\put(10.78,6.8){Initial displ. $x_1^0$ [m]}
\put(10.78,0.4){Initial displ. $x_1^0$ [m]}
\put(0.6,14.45){\rotatebox{90}{Displacement $x_1$ [m]}}
\put(0.6,8.1){\rotatebox{90}{Initial vel. $\dot{x}_1^0$ [m/s]}}
\put(0.6,1.7){\rotatebox{90}{Initial vel. $\dot{x}_1^0$ [m/s]}}
\put(8.5,14.5){\rotatebox{90}{Initial vel. $\dot{x}_1^0$ [m/s]}}
\put(8.5,8.1){\rotatebox{90}{Initial vel. $\dot{x}_1^0$ [m/s]}}
\put(8.5,1.7){\rotatebox{90}{Initial vel. $\dot{x}_1^0$ [m/s]}}
 \put(4.3,18.8){(a)}
 \put(12.2,18.8){(b)}
 \put(4.3,12.4){(c)}
 \put(12.2,12.4){(d)}
 \put(4.3,6.0){(e)}
 \put(12.2,6.0){(f)}
\end{picture}\caption{Basins of attraction of low- (main branch) and high-amplitude (DRC) periodic solutions for $F = 0.15\,$N. (a) Close-up of the frequency response where the frequencies at which the basins of attraction are computed are indicated with dashed lines. (b-f) Basins of attraction for $\omega = 1.67\,$rad/s; $\omega = 1.8\,$rad/s; $\omega = 2\,$rad/s; $\omega = 2.2\,$rad/s, and $\omega = 2.3\,$rad/s, respectively. White and black dots denote the coexisting periodic solutions on the main frequency response and on the DRC, respectively.}
\label{global_f_0_15_isola}
\end{figure}

\begin{figure}[h]
\setlength{\unitlength}{1cm}
\begin{picture}(8,8.5)(0,0)
\put(3,0.6){\includegraphics[width=10truecm]{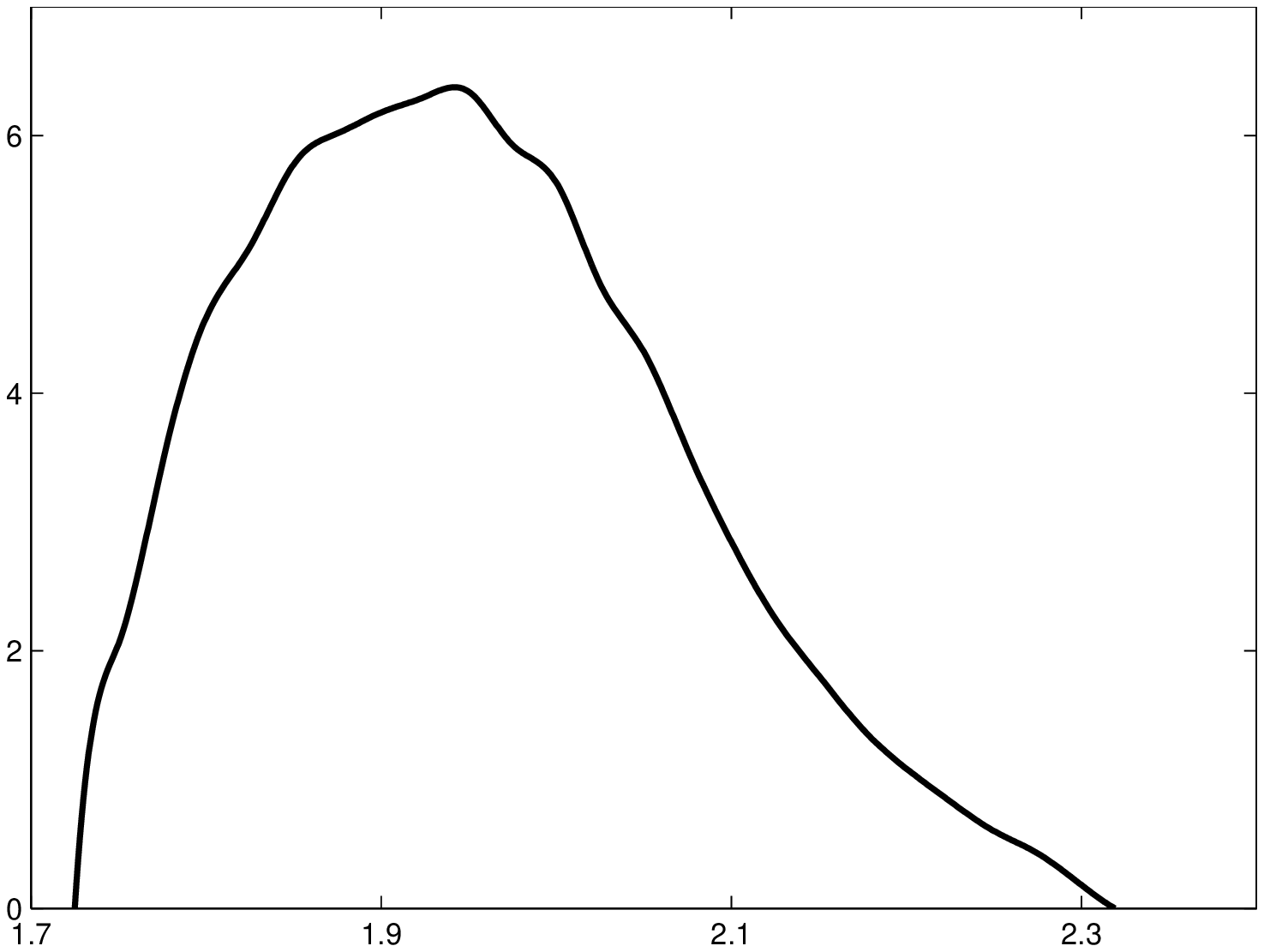}}
\put(6.5,0){Frequency [rad/s]} \put(2.4,2.69){\rotatebox{90}{Basin areas ratio [\%]}}
\end{picture}\caption{Ratio between the areas of the basins of attraction of the DRC and of the main resonance branch for $F = 0.15\,$N.}\label{basin_evol}
\end{figure}

For $F = 0.15\,$N, Figure \ref{global_f_0_15_isola}(a) illustrates that the bistable region lies after the second NS bifurcation on the DRC, i.e., 
in the frequency interval $[1.73-2.32]\,$rad/s. The basins of attraction of the DRC in this interval are found to be very small compared to 
those of the main branch, as shown in Figures \ref{global_f_0_15_isola}(c-f). This finding is confirmed in Figure \ref{basin_evol} where the ratio between the areas of the basins of attraction of the DRC and of the low-amplitude periodic solution does not exceed 6\% for the considered range of initial conditions. 
The basins of attraction of the DRC are also located at a considerable distance from the origin, meaning that high-energy initial conditions 
have to be imparted to the system to excite the DRC. 
Finally, we also verified that the basins of attraction of the DRC remain small for other forcing amplitudes and 
that the low-amplitude solution is the only stable solution in the interval $[1.57-1.73]\,$rad/s, as confirmed in Figure \ref{global_f_0_15_isola}(a).

\subsection{Safe, unsafe and unacceptable NLTVA operations}\label{perfo_reg}

\begin{figure}[h]
\setlength{\unitlength}{1cm}
\begin{picture}(8,9)(0,0)
\put(2.5,0.6){\includegraphics[width=11truecm]{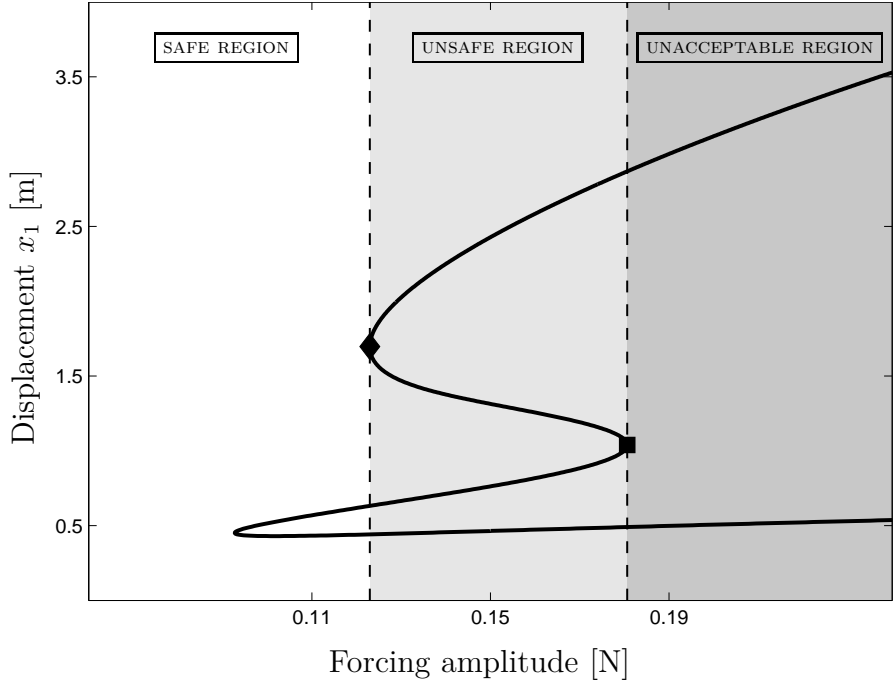}}
\put(6.1,0){Forcing amplitude [N]} \put(1.9,3){\rotatebox{90}{Displacement $x_1$ [m]}}
\put(3.8,8.2){\fbox{\scriptsize{\textsc{safe region}}}}
\put(7.2,8.2){\fbox{\scriptsize{\textsc{unsafe region}}}}
\put(10.15,8.2){\fbox{\scriptsize{\textsc{unacceptable region}}}}
\end{picture}\caption{Performance regions of the NLTVA. The solid line is the projection of branch B of fold bifurcations onto the ($x_1,F$) plane, and the diamond and square markers indicate the apparition and merging of the DRCs, respectively.}\label{LP_perfo_reg}
\end{figure}

Three distinct regions, schematized in Figure \ref{LP_perfo_reg} and characterized based on branch B of fold bifurcations in Figure \ref{LP_tracking_analysis}, are defined for the operation of the NLTVA:
\begin{enumerate}
  \item In the first region, the only branch of periodic solutions is the main branch. Quasiperiodic solutions exist, but they barely degrade NLTVA performance. It is therefore {\it safe} to operate the NLTVA in this region.
  \item The second region presents a large-amplitude DRC. Even if this DRC appears outside the operating frequency range of the NLTVA and if its basins of attraction are small, 
  it is {\it unsafe} to operate the NLTVA in this region.
  \item In the third region, the DRC has merged with the main branch, resulting in a resonance peak of very high amplitude. Even if part of the branch may be unstable due to NS bifurcations, it is {\it unacceptable} to operate the NLTVA in this region.
\end{enumerate}

\section{Sensitivity analysis of the nonlinear tuned vibration absorber}\label{sens_anal}

\subsection{Attenuation performance in the safe region}

The effects of variations of the damping and nonlinear stiffness coefficients on performance in the safe region are now studied. 
Variations of the linear stiffness are not considered herein, because very accurate values of the optimal frequency ratio can easily be 
obtained through small adjustments of the mass ratio (by, e.g., adding small masses on the NLTVA once it is built). 

\begin{figure}[t]
\setlength{\unitlength}{1cm}
\begin{picture}(8,8.5)(0,0)
\put(3,0.6){\includegraphics[width=10truecm]{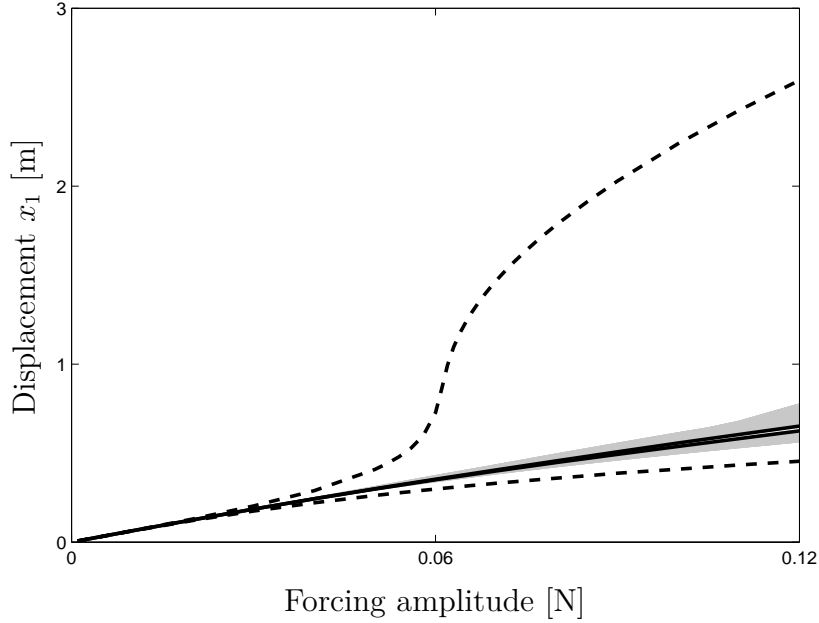}}
\put(6.0,0){Forcing amplitude [N]} \put(2.4,2.6){\rotatebox{90}{Displacement $x_1$ [m]}}
\end{picture}\caption{Performance of the LVTA/NLTVA for increasing forcing amplitudes. The dashed and solid lines depict the amplitude
of the resonances peaks of the Duffing oscillator with an attached LTVA and NLTVA, respectively. The regions in gray show the individual effects of $+15\%$ and $-15\%$ perturbations of the optimal values $\mu_2^{opt}$ and $\beta_3^{opt}$ of the NLTVA.}\label{compa_NLTVA_sensitivity}
\end{figure}

\begin{figure}[!]
\setlength{\unitlength}{1cm}
\begin{picture}(8,6)(0,0)
\put(0.75,0.6){\includegraphics[width=7.0truecm]{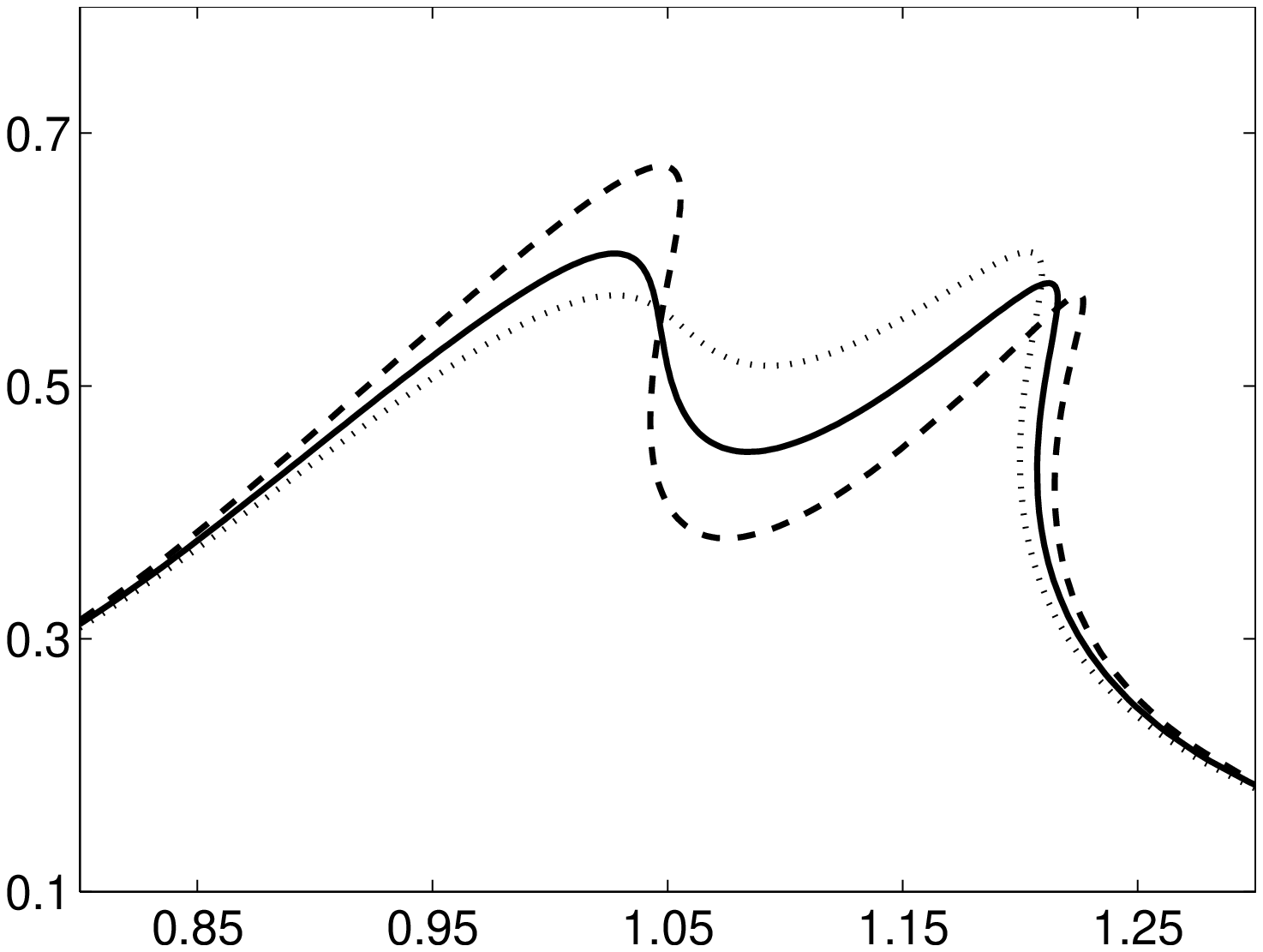}}
\put(8.45,0.6){\includegraphics[width=7.0truecm]{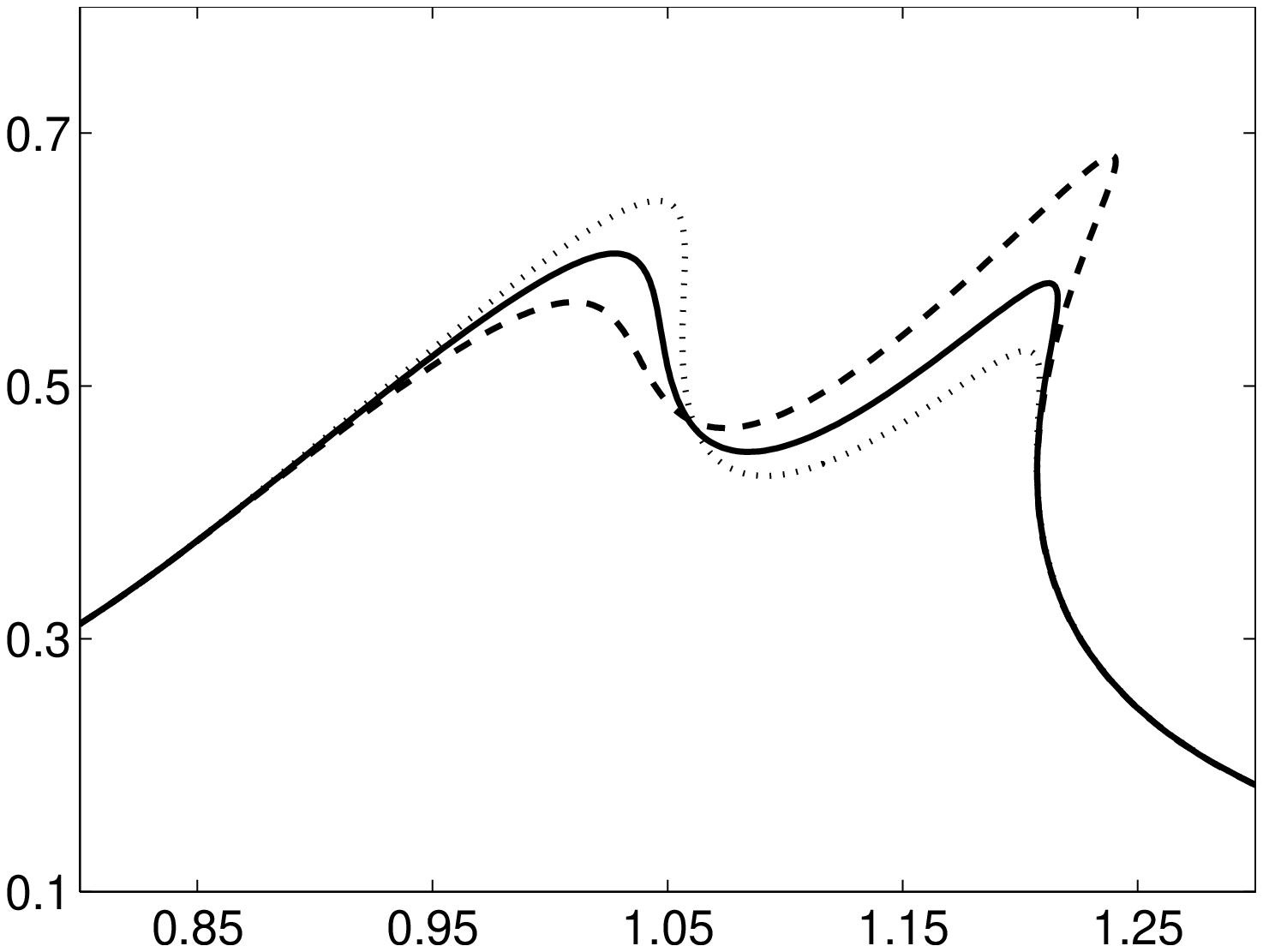}}
\put(6.6,0){Frequency [rad/s]} \put(-0.1,1.4){\rotatebox{90}{Displacement $x_1$ [m]}}
 \put(1.3,5.3){(a)}
 \put(9.1,5.3){(b)}
\end{picture}\caption{Sensitivity of the NLTVA performance with respect to (a) $c_2$ and (b) $k_{nl2}$ for $F = 0.11\,$N. The solid line represents the optimal value, the dashed and dotted lines correspond to variations of -15\% and 15\% with respect to the optimal value, respectively.}\label{mu2_beta3_param_robustness}
\end{figure}

Figure \ref{compa_NLTVA_sensitivity} represents the effects of individual perturbations of $\pm 15\%$ of $c_2$ and $k_{nl2}$ on the amplitude of 
the resonance peaks in the safe region, i.e., until $F = 0.12\,$N. We note that these variations are realistic in view of what was achieved with 
an experimental NLTVA prototype \cite{Chiara}. The NLTVA performance is not significantly degraded and 
remains largely superior to that of the unperturbed LTVA, clearly highlighting the robustness of NLTVA performance in the safe region. For illustration, 
Figure \ref{mu2_beta3_param_robustness} depicts the corresponding frequency responses at $F = 0.11\,$N.

\subsection{Boundaries in NLTVA parameter space}

System (\ref{eomdim2}) is rewritten in dimensionless form. Defining the dimensionless time $\tau=\omega_{n1}t$, where $\omega_{n1}=\sqrt{k_1/m_1}$, applying the transformation $r(t)=x_1(t)-x_2(t)$, and normalizing the system (\ref{eomdim}) using $q_1=x_1/f$ and $q_2=r/f$ (with $f=F/k_1$) yields
\begin{align}\label{eomadim}
\begin{split}
\nonumber
q_1''+2\mu_1q_1'+q_1+\frac{4}{3}\alpha_3q_1^3+2\mu_2\lambda\epsilon q_2'+\lambda^2\epsilon q_2+\frac{4}{3}\epsilon\beta_3 q_2^3={}&\cos\gamma\tau\end{split}\\
\begin{split}
q_2''+2\mu_1q_1'+q_1+\frac{4}{3}\alpha_3q_1^3+2\mu_2\lambda\left(\epsilon+1\right)q_2'+\lambda^2\left(\epsilon+1\right)q_2\\+\frac{4}{3}\left(\epsilon+1\right)\beta_3q_2^3={}&\cos\gamma\tau\end{split}
\end{align}
where prime denotes differentiation with respect to $\tau$, $2\mu_1=c_1/(m_1\omega_{n1})$, $2\mu_2=c_2/(m_2\omega_{n2})$, $\epsilon=m_2/m_1$, $\gamma=\omega/\omega_{n1}$, $\omega_{n2}=\sqrt{k_2/m_2}$, $\lambda=\omega_{n2}/\omega_{n1}$, $\alpha_3 =3k_{nl1}F^2/(4k_1^3)$ and $\beta_3=3k_{nl2}F^2/(4k_1^3\epsilon)$.

Using these dimensionless notations, Equations (\ref{DHrule_dim}-\ref{DHrule_dim2}) can be recast into
\begin{eqnarray}\label{DHrule}
\nonumber
\lambda^{opt}&=&\frac{2}{1+\epsilon}\sqrt{\frac{2\left[16+23\epsilon+9\epsilon^2+2(2+\epsilon)\sqrt{4+3\epsilon} \right]}
{3(64+80\epsilon+27\epsilon^2)}}\\
\nonumber
\mu_2^{opt}&=&\frac{1}{4}\sqrt{\frac{8+9\epsilon-4\sqrt{4+3\epsilon}}{1+\epsilon}}\\
\beta_3^{opt}&=&\frac{2\alpha_3\epsilon}{1+4\epsilon}
\end{eqnarray}
The forcing amplitude $F$ now only appears in the equations through the parameter $\alpha_3$. For instance, for the values in Table \ref{tab:param_sys} and for $F = 0.11$ N, $\alpha_3$ is equal to 0.009075.

\begin{figure}[p]
\setlength{\unitlength}{1cm}
\begin{picture}(8,18)(0,0)
\put(0.60,13.4){\includegraphics[width=7.0truecm]{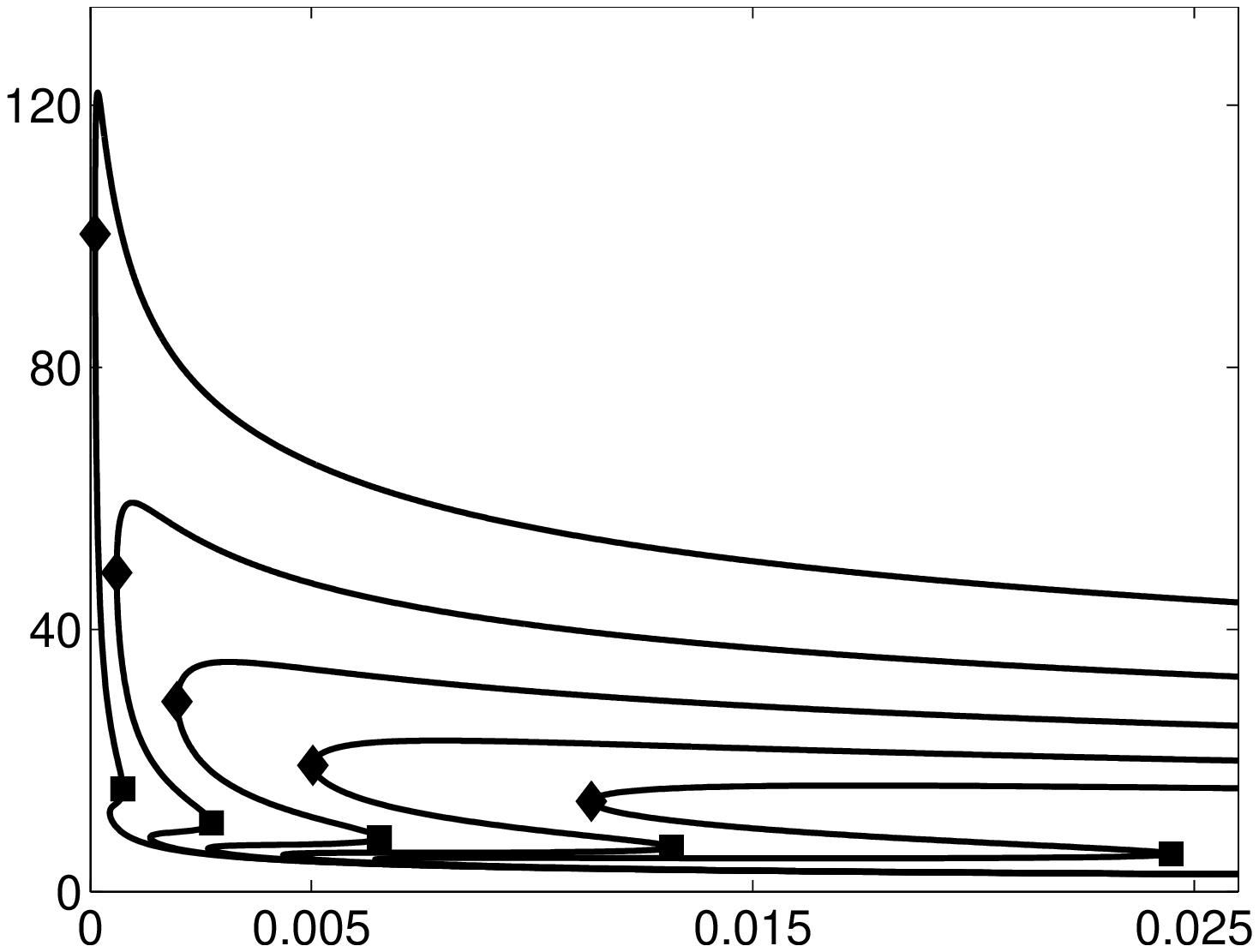}}
\put(1.8,16.5){\scriptsize{$\epsilon = 0.01$}}
\put(1.72,15.7){\scriptsize{$0.02$}}
\put(1.8,15.05){\scriptsize{$0.03$}}
\put(2.5,14.3){\scriptsize{$0.04$}}
\put(5.4,14.08){\scriptsize{$0.05$}}
\put(9.05,13.4){\includegraphics[width=7.0truecm]{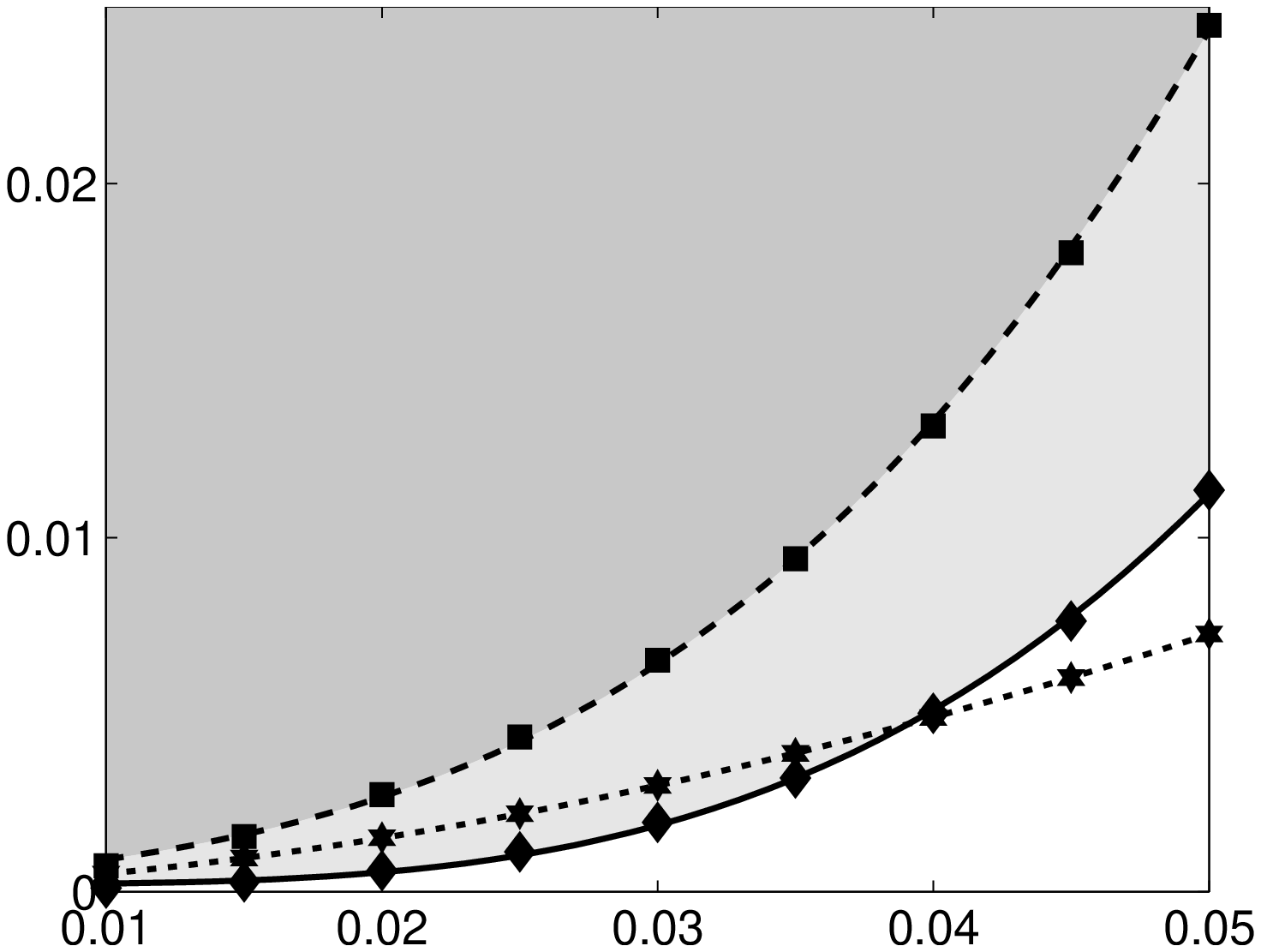}}
\put(14.6,14.2){\fbox{\scriptsize{\textsc{safe}}}}
\put(13.9,15.3){\rotatebox{40}{\fbox{\scriptsize{\textsc{unsafe}}}}}
\put(11.5,16.9){\fbox{\scriptsize{\textsc{unacceptable}}}}
\put(0.60,7){\includegraphics[width=7.0truecm]{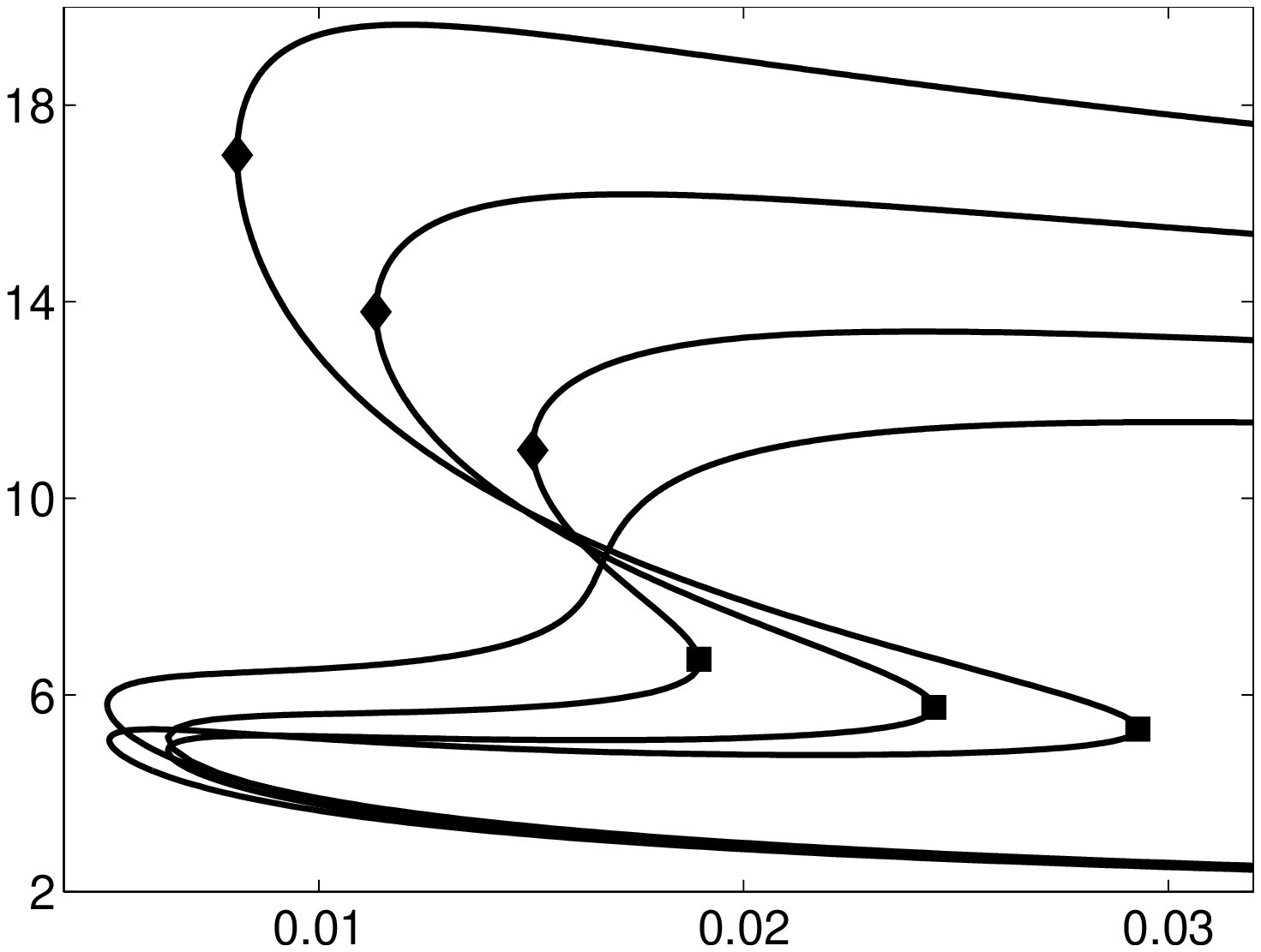}}
\put(2.2,11.7){\scriptsize{$80\%$}}
\put(4,11.3){\scriptsize{$100\%$}}
\put(5.3,10.55){\scriptsize{$125\%$}}
\put(6.0,10.03){\scriptsize{$p_{\mu} = 150\%$}}
\put(9.05,7){\includegraphics[width=7.0truecm]{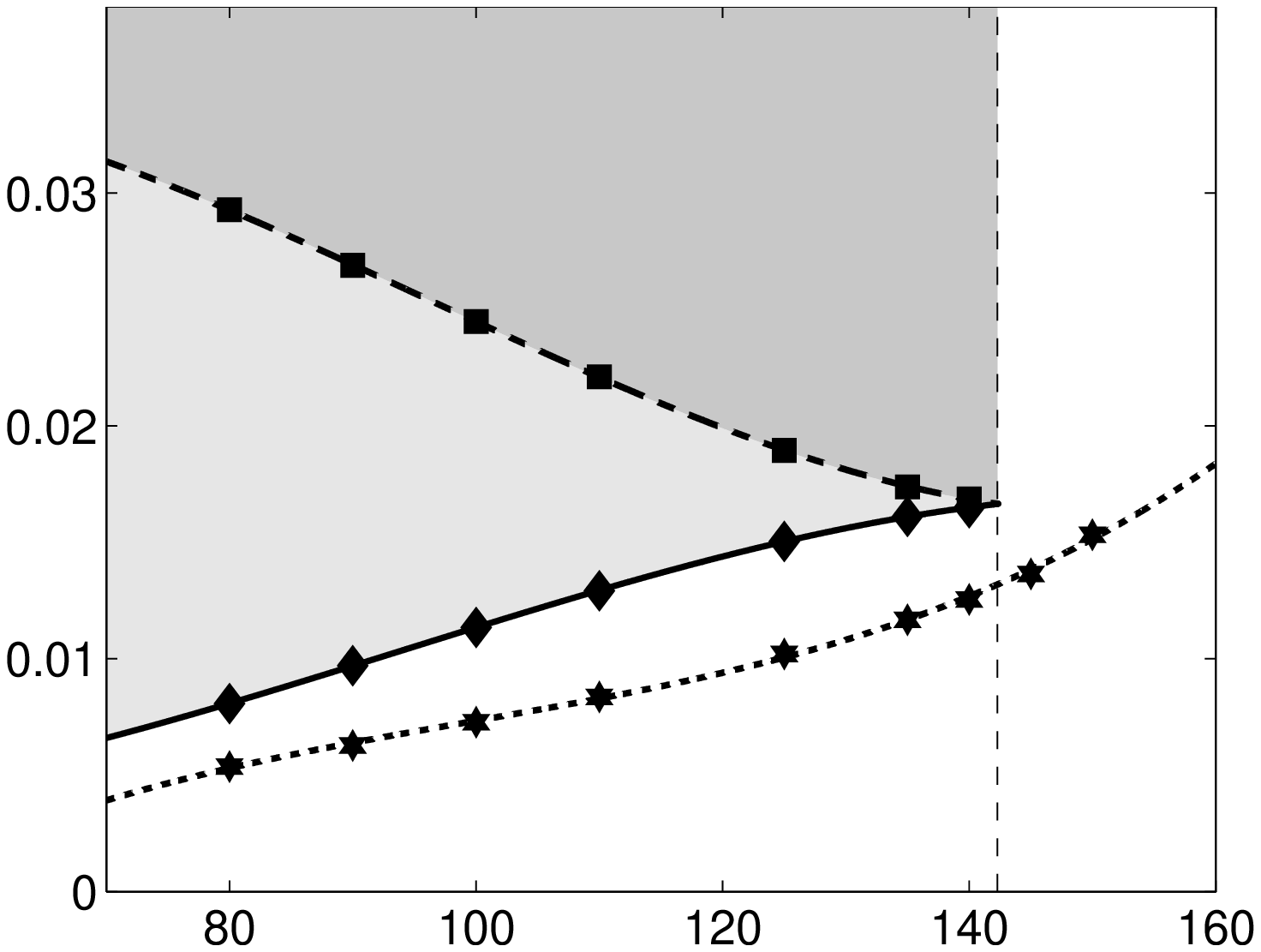}}
\put(13.1,7.9){\fbox{\scriptsize{\textsc{safe}}}}
\put(10.5,9.6){\fbox{\scriptsize{\textsc{unsafe}}}}
\put(12,11.2){\fbox{\scriptsize{\textsc{unacceptable}}}}
\put(0.60,0.6){\includegraphics[width=7.0truecm]{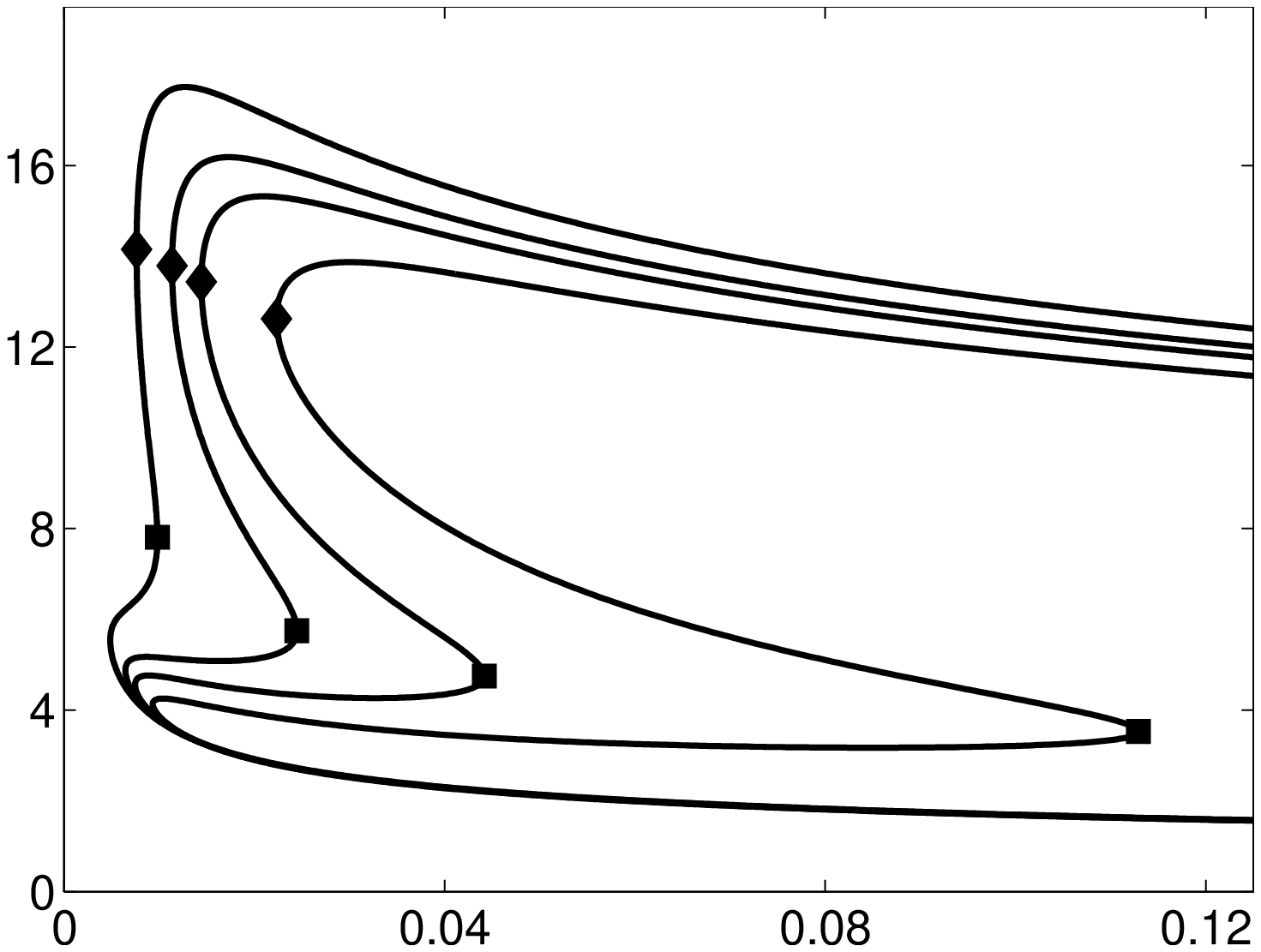}}
\put(1.45,2.4){\scriptsize{$80\%$}}
\put(2.3,2.1){\scriptsize{$100\%$}}
\put(3.4,2.1){\scriptsize{$110\%$}}
\put(5.0,2.4){\scriptsize{$p_{\beta} = 125\%$}}
\put(9.05,0.6){\includegraphics[width=7.0truecm]{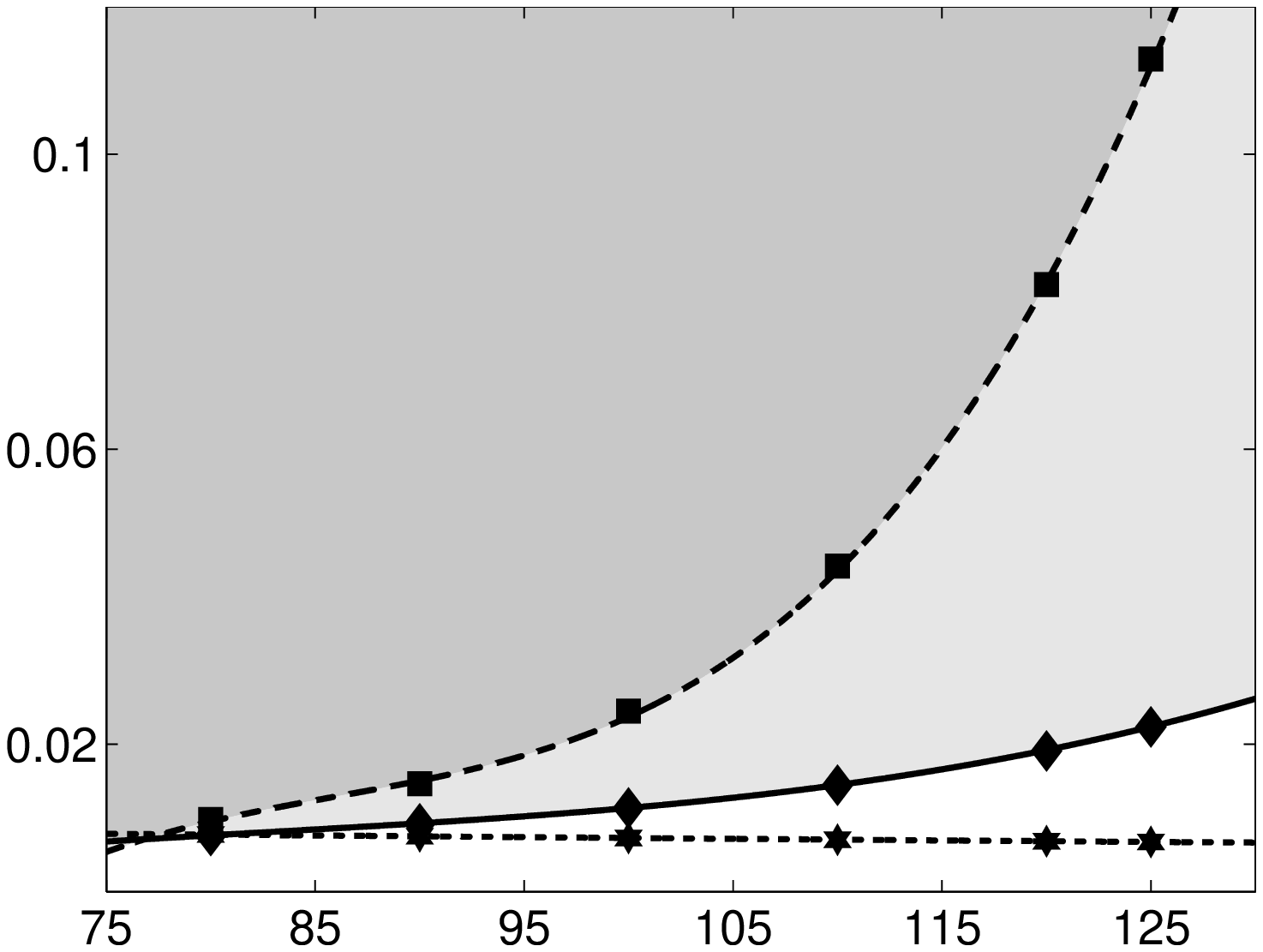}}
\put(15.1,1.42){\fbox{\scriptsize{\textsc{safe}}}}
\put(13.95,2.4){\fbox{\scriptsize{\textsc{unsafe}}}}
\put(11.4,4.2){\fbox{\scriptsize{\textsc{unacceptable}}}}
\put(-0.1,6.8){\rotatebox{90}{Dimensionless displacement $q_1$ [-]}}
\put(8.35,8.3){\rotatebox{90}{Parameter $\alpha_3$ [-]}}
\put(2.85,12.8){Parameter $\alpha_3$ [-]}
\put(2.85,6.4){Parameter $\alpha_3$ [-]}
\put(2.85,0){Parameter $\alpha_3$ [-]}
\put(11.4,12.8){Mass ratio $\epsilon$ [-]}
\put(11.2,6.4){Parameter $p_{\mu}$ [\%]}
\put(11.3,0){Parameter $p_{\beta}$ [\%]}
 \put(1.25,18.2){(a)}
 \put(9.7,18.2){(b)}
 \put(1.00,11.8){(c)}
 \put(9.7,11.8){(d)}
 \put(1.00,5.4){(e)}
 \put(9.7,5.4){(f)}
\end{picture}\caption{Influence of absorber parameters on the regions of NLTVA operation. (a-b) Effect of $\epsilon$; 
(c-d) Effect of $p_{\mu}={\mu_2}/{\mu_2^{opt}}$; (e-f) Effect of $p_{\beta}={\beta_3}/{\beta_3^{opt}}$. First column: projection of the branches of fold bifurcations onto the ($q_1,\alpha_3$) plane. The diamond and square markers indicate the appearance and merging of the DRCs, respectively. Second column: boundaries of the safe, unsafe and unacceptable regions. 
The curves with diamonds, squares and stars represent the boundaries between safe and unsafe regions, unsafe and unacceptable regions, and the onset of quasiperiodic motion, respectively. }\label{fold_reg_summary}
\end{figure}

The influence of parameters $\epsilon$, $\mu_2$ and $\beta_3$ on the boundaries of the different regions of NLTVA operation is now examined. Figure \ref{fold_reg_summary}(a) 
demonstrates the beneficial influence of larger mass ratios. Not only they correspond to resonance peaks of smaller amplitudes (assuming that fold bifurcations occur in the vicinity of the resonance peaks), 
but they also postpone both the appearance and the merging of the DRC to greater values of $\alpha_3$. The corresponding boundaries in 
Figure \ref{fold_reg_summary}(b) illustrate the clear enlargement of the safe region for greater mass ratios. An interesting observation is that, 
for $\epsilon < 4\%$, there are no longer quasiperiodic solutions in the safe region.

Figure \ref{fold_reg_summary}(c) shows that increasing the damping ratio can translate into a fold bifurcation branch that possesses no folding, meaning that the DRC can be completely eliminated. 
Specifically, Figure \ref{fold_reg_summary}(d) confirms that the unsafe and unacceptable regions disappear for ${p}_{\mu}={\mu_2}/{\mu_2^{opt}} \approx 144$\%. Such a detuning of the damping coefficient 
seems interesting, but it is associated with an important decrease in performance in the safe region. 

Figure \ref{fold_reg_summary}(e) illustrates that the nonlinear stiffness coefficient $\beta_3$ has a strong influence on the merging of the DRC. 
Indeed, as plotted in Figure \ref{mu2_beta3_param_robustness}(b), 
greater values of $\beta_3$ reduce the amplitude of the second peak, hence, postponing its merging with the DRC. Increasing $\beta_3$ is also beneficial for delaying the appearance of the DRC and 
enlarging the safe region, as shown in Figure \ref{fold_reg_summary}(f). On the other hand, Figure \ref{mu2_beta3_param_robustness}(b) evidences that a greater $\beta_3$ increases the amplitude of the first resonant peak, 
which, in turn, decreases the NLTVA performance.

\section{Conclusion}\label{conclusion}

In view of the potentially adverse dynamical attractors nonlinear systems can exhibit, the objective of the paper was to identify, and possibly enlarge, 
the safe region of operation of a recently-developed nonlinear absorber, the NLTVA. This was achieved thanks to the combination of several methods 
of nonlinear dynamics, namely the numerical continuation of periodic solutions, bifurcation detection and tracking, and global analysis. Specifically, 
bifurcation tracking proved very useful for determining precisely the creation and elimination of DRCs, which can easily be missed otherwise. It turns out that the best strategy 
to enlarge the safe region while maintaining excellent NLTVA performance is to increase the mass ratio. If it cannot be further increased 
because of practical considerations, an alternative is to increase either damping or the nonlinear coefficient of the absorber.

\section*{Acknowledgments}

The authors T. Detroux, G. Habib, L. Masset and G. Kerschen would like to acknowledge the financial support of the European Union (ERC Starting Grant NoVib
307265).


\begin{thebibliography}{999}

\bibitem{Daraio}
N. Boechler, G. Theocharis, C. Daraio, Bifurcation-based acoustic switching and rectification, {\it Nature Materials} 10, 665-668, 2011.

\bibitem{Antonio}
D. Antonio, D.H. Zanette, D. Lopez, Frequency stabilization in nonlinear mechanical oscillators, {\it Nature Communications} 3, article 806, 2012.

\bibitem{Shaw}
B.S. Strachan, S.W. Shaw, O. Kogan, Subharmonic resonance cascades in a class of coupled resonators, {\it Journal of Computational and Nonlinear Dynamics} 8, Article number 041015, 2013.

\bibitem{Poovarodom}
N. Poovarodom, S. Kanchanosot, P. Warnitchai, Application of non-linear multiple tuned mass dampers to suppress man-induced vibrations of a pedestrian bridge, {\it Earthquake Engineering and Structural Dynamics} 32, 1117-1131, 2003.

\bibitem{Alexander}
N.A. Alexander, F. Schilder, Exploring the performance of a nonlinear tuned mass damper, {\it Journal of Sound and Vibration} 319, 445-462, 2009.

\bibitem{Machado}
M. Febbo, S.P. Machado, Nonlinear dynamic vibration absorbers with a saturation, {\it Journal of Sound and Vibration} 332, 1465-1483, 2013.

\bibitem{Lacarbonara}
N. Carpineto, W. Lacarbonara, F. Vestroni, Hysteretic tuned mass dampers for structural vibration mitigation, {\it Journal of Sound and Vibration} 333, 1302-1318, 2014.

\bibitem{Barton}
D.A.W. Barton, S.G. Burrow, L.R. Clare, Energy harvesting from vibrations with a nonlinear oscillator, {\it Journal of Vibration and Acoustics} 132, 021009, 2010.

\bibitem{Quinn}
D.D. Quinn, A.L. Triplett, A.F. Vakakis, L.A. Bergman, Energy harvesting from impulsive loads using intentional essential nonlinearities, Journal of Vibration and Acoustics 133, 011004, 2011.

\bibitem{Worden}
P.L. Green, K. Worden, K. Atallah, N.D. Sims, The benefits of Duffing-type nonlinearities and electrical optimisation of a mono-stable energy harvester under white Gaussian excitations, Journal of Sound and Vibration 331, 4504-4517, 2012.

\bibitem{Inman}
A. Karami, D.J. Inman, Powering pacemakers from heartbeat vibrations using linear and nonlinear energy harvesting, Applied Physics Letters 100, 042901, 2012.

\bibitem{Book}
A.F. Vakakis, O. Gendelman, L.A. Bergman, D.M. McFarland, G. Kerschen, Y.S. Lee, {\it Nonlinear Targeted Energy Transfer in Mechanical and Structural Systems}, {\it Springer}, Series: Solid Mechanics and Its Applications, 2009.

\bibitem{NLD} G. Kerschen, J.J. Kowtko, D.M. McFarland, L.A. Bergman, A.F. Vakakis, Theoretical and experimental study of multimodal tet in a system of
coupled oscillators, {\it Nonlinear Dynamics} 47, 285-309, 2007.

\bibitem{SIAM} G. Kerschen, Y.S. Lee, A.F. Vakakis, D.M. McFarland, L.A. Bergman, Irreversible passive energy transfer in coupled oscillators with
essential nonlinearity, {\it SIAM Journal on Applied Mathematics} 66, 648-679, 2006.

\bibitem{Nucera}
F. Nucera, A.F. Vakakis, D.M. McFarland, L.A. Bergman, G. Kerschen, Targeted energy transfers in vibro-impact oscillators for seismic mitigation, {\it Nonlinear Dynamics} 50 (2007), 651-677.

\bibitem{Hubbard}
S.A. Hubbard, D.M. McFarland, L.A. Bergman, A.F. Vakakis, Targeted energy transfer between a model flexible wing and nonlinear energy sink, {\it Journal of Aircraft} 47, 1918-1931, 2010.

\bibitem{Lamarque}
B. Vaurigaud, L.I. Manevitch, C.H. Lamarque, Passive control of aeroelastic instability in a long span bridge model prone to coupled flutter using targeted energy transfer, {\it Journal of Sound and Vibration} 330, 2580-2595, 2011.

\bibitem{Cochelin}
R. Bellet, B. Cochelin, P. Herzog, P.O. Mattei, Experimental study of targeted energy transfer from an acoustic system to a nonlinear membrane absorber,
{\it Journal of Sound and Vibration} 329, 2768-2791, 2010.

\bibitem{Gourc}
E. Gourc, S. Seguy, G. Michon, A. Berlioz, Chatter control in turning process with a nonlinear energy sink,
{\it Advanced Materials Research} 698, 89-98, 2013.

\bibitem{NLTVA1}
G. Habib, T. Detroux, R. Viguié, G. Kerschen, Nonlinear generalization of Den Hartog's equal-peak method,
{\it Mechanical Systems and Signal Processing} 52-53, 17-28, 2015.

\bibitem{ENOC}
G. Habib, G. Kerschen, Stability and bifurcation analysis of a Van der Pol-Duffing oscillator with a nonlinear tuned vibration absorber,
European Nonlinear Dynamics Conference, Vienna, Austria, 2014, available online at {\it http://orbi.ulg.ac.be/handle/2268/167407}.

\bibitem{Staros}
Y. Starosvetsky, O.V. Gendelman, Response regimes of linear oscillator coupled
to nonlinear energy sink with harmonic forcing and frequency detuning, {\it Journal of Sound and Vibration}
315 (2008), 746-765.

\bibitem{Shaw4}
J. Shaw, S.W. Shaw, A.G. Haddow, On the Response of the Non-linear Vibration Absorber, {\it International Journal of Non-linear Mechanics} 24, 281-293, 1989.


\bibitem{staros2}
Y. Starosvetsky, O.V. Gendelman, Vibration absorption in systems with a nonlinear energy sink: Nonlinear damping, {\it Journal of Sound and Vibration} 324,
916-939, 2009.

\bibitem{Gourc2}
E. Gourc, G. Michon, S. Seguy, A. Berlioz, Experimental investigation and design optimization of targeted energy transfer under periodic forcing,
{\it Journal of Vibration and Acoustics} 136, 021021, 2014.

\bibitem{Singh}
C. Duan, R. Singh, Isolated sub-harmonic resonance branch in the frequency response of an oscillator with slight asymmetry in the clearance, {\it Journal of Sound and Vibration} 314, 12-18, 2008.

\bibitem{Stepan}
D. Takacs, G. Stepan, S.J. Hogan, Isolated large amplitude periodic motions of towed rigid wheels, {\it Nonlinear Dynamics} 52, 27-34, 2008.

\bibitem{Grolet}
A. Grolet, E. Sarrouy, F. Thouverez, Global and bifurcation analysis of a structure with cyclic symmetry, {\it International Journal of Non-linear Mechanics} 46, 727-737, 2011.

\bibitem{Asami} T. Asami, O. Nishihara, Closed-form exact solution to H$\infty$ optimization of dynamic vibration absorbers
(application to different transfer functions and damping systems), {\it Journal of Vibration and Acoustics} 125, 398-405, 2003.

\bibitem{DenHartog} J.P. Den Hartog, {\it Mechanical Vibrations}, McGraw-Hill, New York, 1934.

\bibitem{Brock} J.E. Brock, A note on the damped vibration absorber, {\it Journal of Applied Mechanics} 13, pp. A284, 1946.

\bibitem{Orlando}
T. Detroux, L. Renson, L. Masset, G. Kerschen, The harmonic balance method for bifurcation analysis nonlinear mechanical systems,
International Modal Analysis Conference, Orlando, Florida, USA, 2015, available at {\it http://orbi.ulg.ac.be/handle/2268/173556}.



\bibitem{Afrai}
V.S. Afraimovich, L.P. Shilnikov, Invariant two-dimensional tori, their breakdown and stochasticity, {\it Transactions of the American Mathematical Society} 149, 201-212, 1991.

\bibitem{Matcont}
A. Dhooge, W. Govaerts, Y. U. Kuznetsov, MATCONT: a MATLAB package for numerical bifurcation analysis of ODEs, {\it ACM Transactions on Mathematical Software} 29, 141-164, 2003.

\bibitem{Dick}
R.P. Eason, A.J. Dick, S. Nagaragaiah, Numerical investigation of coexisting high and low amplitude responses and safe basin erosion for a coupled linear oscillator and nonlinear absorber system, {\it Journal of Sound and Vibration} 333, 3490-3504, 2014.

\bibitem{Chiara}
C. Grappasonni, G. Habib, T. Detroux, G. Kerschen, Experimental demonstration of a 3D-printed nonlinear tuned vibration absorber, {\it 33rd International Modal Analysis Conference}, Orlando, USA, 2015, available at {\it http://orbi.ulg.ac.be//handle/2268/173803}.




%

%
%

%



\end{thebibliography}
\end{document}